\documentclass{amsart}
\usepackage{setspace}

\usepackage{amsmath}
\usepackage{amsthm}
\usepackage{amssymb}
\usepackage{color}
\usepackage[latin1]{inputenc} 
\usepackage{graphicx}

\newtheorem{theorem}{Theorem}[section]
\newtheorem{lemma}[theorem]{Lemma}
\newtheorem{proposition}[theorem]{Proposition}
\newtheorem{corollary}[theorem]{Corollary}
\newtheorem{remark}[theorem]{Remark}

\newtheorem{definition}[theorem]{Definition}

\numberwithin{equation}{section}

\begin{document}
\title[Entropy and divergence in number theory]
{Entropy and divergence in number theory}

\author{Daniel C. Mayer}
\address{Karl-Franzens University Graz\\
 Graz, Austria}
\email{algabraic.number.theory@algebra.at; quantum.algebra@icloud.com}

\author{Nicu\c{s}or Minculete}
\address{Faculty of Mathematics and Computer Science, Transilvania University\\
 Iuliu Maniu street 50, Bra\c{s}ov 500091, Romania}
\email{minculete.nicusor@unitbv.ro; minculeten@yahoo.com}

\author{Diana Savin}
\address{Faculty of Mathematics and Computer Science, Transilvania University\\
 Iuliu Maniu street 50, Bra\c{s}ov 500091, Romania}
\email{diana.savin@unitbv.ro; dianet72@yahoo.com}

\author{Vlad Monescu}
\address{Faculty of Mathematics and Computer Science, Transilvania University\\
 Iuliu Maniu street 50, Bra\c{s}ov 500091, Romania}
\email{monescu@unitbv.ro}

\subjclass[2010]{Primary: 28D20, 11A51, 11A25, 11R37, 11R29, 11R11; Secondary: 11S15, 20D15, 47B06, 94A17}
\keywords{entropy, divergence, numbers, ideals, ramification theory in algebraic number fields, \(3\)-class field towers, imaginary quadratic fields, Schur \(\sigma\)-groups.}

\date{16 September 2025}

\begin{abstract}
We obtain inequalities involving the entropy of a positive integer and the divergence of two positive integers,
respectively the entropy of an ideal and the divergence of two ideals in a ring of algebraic integers.
Among the important results,
we show that the minimal entropy arises for sharp localization,
and the maximal entropy occurs for equidistribution.
We also study other interesting estimates of entropy and divergence for numbers and for ideals.
Finally, we determine the entropies of probability distributions on infinite trees of Schur \(\sigma\)-groups,
which are realized by \(3\)-class field tower groups of imaginary quadratic number fields.
\end{abstract}

\maketitle

\section{Introduction and Preliminaries}
\label{s:Intro}

\noindent 
Let \(\mathcal{G}=(V,E)\) be a countable digraph with 
vertices \(v\in V\) and directed edges \(e=(v_0\to v_1)\in E\subset V\times V\).
Denote by \((0;1\rbrack\) the interval of real numbers \(0<r\le 1\), and
suppose \(p:\,S\to(0;1\rbrack\), \(v\mapsto p(v)\),
is a probability distribution with \(\sum_{v\in S}\,p(v)=1\)
on some subset \(S\subset V\).
By means of the natural logarithm \(\log\),
the \textit{entropy} of \(p\) is then defined to be the weighted sum
\(H(p)=\sum_{v\in S}\,p(v)\cdot\log(p(v))\),
provided it converges.
This general setting will be applied in two arithmetic situations.

Firstly, in algebraic number theory,
to the lattice of integral ideals, viewed as a digraph \(\mathcal{L}=(\mathcal{I}_K,\mathcal{D})\)
with integral ideals \(I\in\mathcal{I}_K\) of an algebraic number field \(K\) as vertices, and
divisor relations \((I_0\mid I_1)\in\mathcal{D}\subset\mathcal{I}_K\times\mathcal{I}_K\) as directed edges.

Secondly, in class field theory,
to a descendant tree of finite \(3\)-groups, viewed as a digraph \(\mathcal{T}=(V,E)\)
containing Galois groups \(G=\mathrm{Gal}(K^\infty/K)\in V\) of maximal unramified \(3\)-extensions \(K^\infty\)
of imaginary quadratic fields \(K\) as vertices, and
descendant-ancestor relations \((G_0\to G_1)\in E\subset V\times V\)
by means of the natural projection \(G_0\to G_0/\gamma_c(G_0)\) to the last non-trivial lower central quotient,
where \((\gamma_i(G_0))_{1\le i\le c+1}\) denotes the lower central series and
\(c\) the nilpotency class of \(G_0\).

The layout of the article is as follows.
In \S\ 
\ref{s:RationalNumberField},
we take the rational number field \(K=\mathbb{Q}\) as base field,
and we identify integral ideals of \(\mathbb{Q}\) with positive integers \(\mathbb{N}\),
in particular, we put \(S\subset\mathbb{P}\) the finite subset of prime numbers
dividing a fixed positive integer \(n\in\mathbb{N}\).
In \S\ 
\ref{s:AlgebraicNumberField},
we investigate the finite subset \(S\subset\mathbb{P}_K\) of prime ideals
dividing a fixed integral ideal \(I\in\mathcal{I}_K\)
in a general number field \(K\).
Finally, in \S\
\ref{s:ClassFieldTheory},
we study the distribution of Schur \(\sigma\)-groups,
arising as Galois groups \(G=\mathrm{Gal}(K^\infty/K)\) of \(3\)-class field towers \(K^\infty\)
of imaginary quadratic fields \(K=\mathbb{Q}(\sqrt{d})\), \(d<0\),
with elementary bicyclic \(3\)-class group \(\mathrm{Cl}_3(K)=(\mathbb{Z}/3\mathbb{Z})^2\),
in infinite subsets \(S\subset V\),
by means of probability measures introduced by Boston, Bush and Hajir in
\cite{BBH2017,BBH2021}.
Here, the entropy is an infinite series with dominated convergence by a geometric series.

In the cases with finite set \(S\),
viewed as a subset of the infinite set \(\mathbb{P}_K\) of all non-archimedean places of \(K\),
minimal entropy \(H(p)=0\) characterizes the maximal amount of order with sharp localization in a singleton set \(S\),
and maximal entropy \(H(p)=\log(\lvert S\rvert)\) is associated with equidistribution,
as a description of maximal disorder.
Our perspective of entropy complements the following well-known classical viewpoints.

In information theory, the entropy is defined as a measure of uncertainty. Over the years, various authors have introduced several types of entropies. One of the most well-known types of entropy is Shannon's entropy $H_S$. This has been defined
 for a probability distribution ${\bf p}=\{p_1,...,p_r\}$ in the following way
$$H_S({\bf p})= - \sum_{i=1}^{r}p_i\cdot \log p_i,$$
where $\sum_{i=1}^{r}p_i=1$ and $0<p_i\le 1$ for all $i=1,\ldots,r$.\\
The most important properties of Shannon's entropy are:
\begin{enumerate} 
 \item
$H_S({\bf pq})=H_S({\bf p})+H_S({\bf q})$, where ${\bf p}=\{p_1,...,p_r\}$, ${\bf q}=\{q_1,...,q_r\}$ and ${\bf pq}=\{p_1q_1,...,p_1q_r,...,p_rq_1,...,p_rq_r\}$ (the additivity);
 \item
$H_S(p_1,p_2,...,p_r)=H_S(p_1+p_2,p_3,...,p_r)+(p_1+p_2)H_S(\frac{p_1}{p_1+p_2},\frac{p_2}{p_1+p_2})$\\
(the recursivity).
 \end{enumerate}
In physics, the entropy has many physical implications as the amount of \lq\lq disorder\rq\rq\ of a system. Entropy is useful in characterizing the behavior of stochastic processes because it represents the uncertainty and disorder of the process. In \cite{DST}, De Gregorio, S\' anchez and Toral defined the block entropy (based on Shannon entropy), which can determine the memory for modeled systems as Markov chains of arbitrary finite order.

Cover and Thomas \cite{CT} introduced the relative entropy (or Kullback--Leibler distance) between two probability distributions ${\bf p}=\{p_1,...,p_r\}$  and ${\bf q}=\{q_1,...,q_r\}$ as follows:
$$
D({\bf p}||{\bf q}):=-\sum_{i=1}^{r}p_i\cdot \log \frac{q_i}{p_i}=\sum_{i=1}^{r}p_i\cdot \log \frac{p_i}{q_i},
$$ 
where $\sum_{i=1}^{r}p_i=1$ and $\sum_{i=1}^{r}q_i=1$ and $0<p_i,q_i\le 1$ for all $i=1,\ldots,r$.

Peculiarities of number theory, related to the factorization of an integer, are given by Dujella in \cite{Duj}. Let $n$ be a positive integer, $n\geq 2$. Minculete and Pozna \cite{Minculete1} introduced the notion of entropy of  $n$  as follows: if  $n=p^{\alpha_{1}}_{1}p^{\alpha_{2}}_{2}\cdots p^{\alpha_{r}}_{r}$, where $r,\alpha_1,\alpha_2,\ldots,\alpha_r\in\mathbb{N}^\ast$ and $p_1,p_2,\ldots,p_r$ are distinct prime positive integers (this representation of $n$ is unique, according to the Fundamental Theorem of Arithmetic), then the entropy of $n$ is:
\begin{equation}\label{1.1}
H\left(n\right)= - \sum_{i=1}^{r}p\left(\alpha_{i}\right)\cdot \log \: p\left(\alpha_{i}\right),    
\end{equation}
where $\log$ denotes the natural logarithm and $p(\alpha_i)=\frac{\alpha_i}{\Omega(n)}$ is a particular probability distribution associated to $n$. By convention, $H(1)=0$.\\
An equivalent form of the entropy of $n\geq 2$ was introduced in \cite{Minculete1} as follows:
\begin{equation}\label{1.2}
H\left(n\right)=\log \: \Omega\left(n\right) - \frac{1}{\Omega\left(n\right)}\cdot \sum_{i=1}^{r}\alpha_{i}\cdot \log \: \alpha_{i},  
\end{equation}
where $\Omega\left(n\right)=\alpha_{1}+\alpha_{2}+...+\alpha_{r}$.

Let $n$ be a positive integer, $n\geq 2$. We denote by $\omega(n)$ the number of distinct prime factors of $n$. In \cite{Minculete1}, the authors defined the Kullback--Leibler distance between two positive integer numbers $n,m\geq 2$ with  factorizations $n=p^{\alpha_{1}}_{1}p^{\alpha_{2}}_{2}\cdots p^{\alpha_{r}}_{r}$ and $m=q^{\beta_{1}}_{1}q^{\beta_{2}}_{2}\cdots q^{\beta_{r}}_{r}$, where the prime factors are arranged in ascending order and $\omega(n)=\omega(m)$, as follows:                                          
\begin{equation}\label{1.3}
D(n||m):=-\sum_{i=1}^{r}p(\alpha_i)\cdot \log \frac{p(\beta_i)}{p(\alpha_i)}, 
\end{equation}
where $p(\alpha_i)=\frac{\alpha_i}{\Omega(n)}$  and  $p(\beta_i)=\frac{\beta_i}{\Omega(m)}$, for every $i\in\{1,2,\ldots,r\}$.
It is clear that  $\sum_{i=1}^{r}p(\alpha_i)=1$ and $\sum_{i=1}^{r}p(\beta_i)=1$ are probability distributions.\\
Formula \eqref{1.3} is equivalent to

\begin{equation}\label{1.4}
D(n||m)=\log\frac{\Omega(m)}{\Omega(n)}-\frac{1}{\Omega(n)}\sum_{i=1}^{r}\alpha_i\cdot \log \frac{\beta_i}{\alpha_i}. 
\end{equation}

\noindent
In \cite{Minculete1}, the authors found crucial properties of the entropy of a positive integer.

\begin{proposition}
\label{onedotone}
The following statements hold generally:
\begin{enumerate} 
 \item
$0\leq H\left(n\right)\leq \log \: \omega\left(n\right)$, for all $n\in \mathbb{N},$\;$n\geq 2;$  
\item
If $n =p^{\alpha}$, with $\alpha$ a positive integer and $p$ a positive prime integer, then $H\left(n\right)=0$
(minimal entropy and high order for sharp localization);
\item
If $n =p_{1}\cdot p_{2}\cdot\ldots\cdot p_{r},$ with $p_{1}, p_{2},\ldots, p_{r}$ distinct positive prime integers, then $H\left(n\right)=\log\:\omega\left(n\right)$
(maximal entropy for equidistribution);
\item
If $n =\left(p_{1}\cdot p_{2}\cdot\ldots\cdot p_{r}\right)^{\alpha},$ with $\alpha$ a positive integer and $p_{1}, p_{2},\ldots, p_{r}$ distinct positive prime integers, then also $H\left(n\right)=log\:\omega\left(n\right)$ (disorder).
\end{enumerate}
\end{proposition}

In \cite{MS1}, Minculete and Savin obtained the following properties involving the entropy and divergence of positive integers.

\begin{proposition}
\label{onedottwo}
Let $n$ and $m$ be two positive integers,  $n,m\geq 2$. Then the following statements are true:
\begin{enumerate} 
 \item
If $n=m$, then we have $D(n||m)=0$;
 \item
If the unique factorizations (in a product of prime factors) of $n$ and $m$ are  $n=p^{\alpha_{1}}_{1}p^{\alpha_{2}}_{2}\ldots p^{\alpha_{r}}_{r}$ and $m=q^{\alpha_{1}}_{1}q^{\alpha_{2}}_{2}\ldots q^{\alpha_{r}}_{r}$, then $D(n||m)=D(m||n)=0$;
 \item
In general, however, $D(n||m)\neq D(m||n)$;
 \item
$H(n^{\alpha})=H(n)$, for any positive integer $\alpha$;
\item
If $\omega(m)=\omega(n)$, then
$D(n||m)=H(m)-H(n)+\sum_{i=1}^{r}\left(\frac{\beta_i}{\Omega(m)}-\frac{\alpha_i}{\Omega(n)}\right)\log \beta_i$.
\end{enumerate}
\end{proposition}

Let $K$ be an algebraic number field. Its ring of algebraic integers is denoted by  $\mathcal{O}_{K}$.  Let $I\neq(0)$ be an ideal of $\mathcal{O}_{K}$.  According to the fundamental theorem of Dedekind rings, $I\neq(1)$ is represented uniquely in the form $I=P^{e_{1}}_{1}\cdot P^{e_{2}}_{2}\cdot\ldots\cdot P^{e_{g}}_{g}$,
where $P_{1},P_{2},\ldots,P_{g}$ are distinct prime ideals of the ring $\mathcal{O}_{K}$ and  $e_{1},e_{2},\ldots,e_{g}$ are positive integers. Let  $ \Omega(I)=e_{1}+e_{2}+\ldots+e_{g}$. Note that \(\mathcal{O}_{K}=(1)$ and $\Omega(\mathcal{O}_{K})=0$.\\
Minculete and Savin \cite{MS} introduced the following notion of entropy of an ideal of the ring $\mathcal{O}_{K}$:
\begin{definition}
\label{onedotthree}
(Definition 1 from \cite{MS}).  Let $I\neq \left(1\right)$ be an ideal of the ring  $\mathcal{O}_{K}$, decomposed as above. We define the entropy of the ideal $I$ in the following way:
$$ H\left(I\right):= - \sum_{i=1}^{g}\frac{e_{i}}{\Omega(I)} \log \: \frac{e_{i}}{\Omega(I)}. $$
\end{definition}

\noindent
In \cite{MS}, the authors also gave an equivalent form of the entropy of the ideal $I\neq(1)$:
\begin{equation}
\label{1.5}
H\left(I\right)=\log \: \Omega\left(I\right) - \frac{1}{\Omega\left(I\right)}\cdot \sum_{i=1}^{g}e_{i}\cdot \log \: e_{i}.  \tag{1.5}
\end{equation}

\noindent
Minculete and Savin \cite{MS1} introduced the notion of the divergence of two ideals of the ring $\mathcal{O}_{K}$ as follows:
\begin{definition}
\label{onedotfour}
(Definition 3.2 from \cite{MS1}).  
Let $I,J\neq \left(1\right)$ be two ideals of the ring  $\mathcal{O}_{K}$, uniquely decomposed as $I=P^{e_{1}}_{1}\cdot P^{e_{2}}_{2}\cdot\ldots\cdot P^{e_{g}}_{g}$ and $J=Q^{f_{1}}_{1}\cdot Q^{f_{2}}_{2}\cdot\ldots\cdot Q^{f_{g}}_{g}$, with $e_{1}, e_{2},\ldots,e_{g}, f_{1}, f_{2},\ldots, f_{g}, $ positive integers, $P_{1}, P_{2},\ldots,P_{g}$ distinct prime ideals of the ring $\mathcal{O}_{K}$ and $Q_{1}, Q_{2},\ldots,Q_{g}$ distinct prime ideals of the ring $\mathcal{O}_{K}.$  Let  $ \Omega\left(I\right) =  e_{1}+ e_{2}+\ldots+ e_{g}$ and $ \Omega\left(J\right) =  f_{1}+ f_{2}+\ldots+ f_{g}.$
We define the divergence of the ideals $I$ and $J$ in the following manner:
\begin{equation}
D\left(I||J\right):= \log\frac{\Omega(J)}{\Omega(I)}-\frac{1}{\Omega(I)}\sum_{i=1}^{g}e_i\cdot \log \frac{f_i}{e_i},    \tag{1.6}
\end{equation}
where $e_i\leq e_j$ and $f_i\leq f_j$ when $i<j$, $i,j\in\{1,\ldots,g\}$.
\end{definition}

The extension of some properties of the natural numbers to ideals was recently given in \cite{MS_Ex}, using the exponential divisors of a natural number and the exponential divisors of an ideal. 

In this article we obtain certain inequalities involving the entropy of a positive integer and divergence of two positive integers, respectively the entropy of an ideal and divergence of two ideals of a ring of algebraic integers. In section 2 we present some inequalities related to the entropy of a positive integer or the divergence of two positive integers. One of the important results shows that an integer number $n\geq 2$ has zero entropy if and only if $n=p^{\alpha}$, where $\alpha$ is a positive integer and $p$ is a prime number. In section 3 we present some inequalities involving the entropy of an ideal of a ring of algebraic integers or the divergence of two ideals of a ring of algebraic integers. Among the important results we studied that if $K$ is an algebraic number field with $\mathcal{O}_{K}$ its ring of algebraic integers and $J\ne (1)$ is an ideal of the ring $\mathcal{O}_{K},$ then $H\left(J\right)=0$ if and only if $J=P^{\alpha}$, where $P$ is a prime ideal of the ring  $\mathcal{O}_{K}$ and $\alpha$ is a positive integer. We also studied how entropy or divergence changes with different factorizations of numbers and how entropy or divergence changes with different factorizations of ideals.

\section{Some inequalities related to the entropy of a positive integer and to the divergence of two positive integers}
\label{s:RationalNumberField}

\noindent 
Let \(K=\mathbb{Q}\),
and consider the divisor lattice \(\mathbb{N}\).
To begin with, we highlight a fundamental property of the entropy \(H(n)\) of a positive integer \(n\in\mathbb{N}\).
\begin{proposition}
\label{twodotone_1.0}
Let $n\geq 2$ be an integer number. Then $H(n)=0$ if and only if $n=p^{\alpha}$, where $\alpha$ is a positive integer and $p$ is a prime number.
\end{proposition}
\begin{proof}
According to the Fundamental Theorem of Arithmetic, 
an integer $n\ge 2$ has a unique representation $n=p_1^{\alpha_1}p_2^{\alpha_2}\cdots p_r^{\alpha_r}$
with at least one prime factor, that is, $r\ge 1$,
distinct prime divisors $p_1<p_2<\ldots<p_r$ arranged in ascending order,
and non-zero exponents $\alpha_i\ge 1$ for $i=1,\ldots,r$.
By Formula \eqref{1.1}, the entropy of $n$ is defined as
$H(n)=-\sum_{i=1}^{r}\,p(\alpha_i)\cdot\log \: p(\alpha_i)$,
where $\log$ is the natural logarithm, and the $p(\alpha_i)=\frac{\alpha_i}{\Omega(n)}$
with $\Omega(n)=\sum_{i=1}^{r}\alpha_i$
form a particular probability distribution associated to $n$. 
Minimal entropy is equivalent to sharp localization:

Sufficiency ($\Longleftarrow$):
If $n=p^\alpha$, then $r=1$, $p=p_1$, $\alpha=\alpha_1$, $\Omega(n)=\alpha$, and $p(\alpha)=\frac{\alpha}{\Omega(n)}=1$,
whence $H(n)=-p(\alpha)\cdot\log \: p(\alpha)=-1\cdot\log \: 1=0$.

Necessity ($\Longrightarrow$):
Suppose that $H(n)=0$ for $n=p_1^{\alpha_1}p_2^{\alpha_2}\cdots p_r^{\alpha_r}$.
If we had more than one prime factor, that is, $r\ge 2$,
then $\Omega(n)=\alpha_1+\alpha_2+\ldots+\alpha_r>\alpha_i$, $p(\alpha_i)=\frac{\alpha_i}{\Omega(n)}<1$, and $\log \: p(\alpha_i)<0$,
for each $i=1,\ldots,r$.
Consequently, the entropy $H(n)=-\sum_{i=1}^{r}\,p(\alpha_i)\cdot\log \: p(\alpha_i)$
would be a sum of at least two positive terms $p(\alpha_i)\cdot(-\log \: p(\alpha_i))>0$,
in contradiction to the assumption that $H(n)=0$.
Thus $r=1$ and $n=p^\alpha$ with $p=p_1$, $\alpha=\alpha_1$.
\end{proof}

\noindent
We consider the natural number $n=p^{\alpha_{1}}_{1}p^{\alpha_{2}}_{2}...p^{\alpha_{r}}_{r}>1$.
We want to study the entropy when $\alpha_i\in\{1,2\}$ for all $i\in\{1,...,r\}$, i.e., for a number $n=p^{2}_{1}p^{2}_{2}...p_s^2p_{s+1}...p_{r}>1$, with $1\leq s\leq r$.
Therefore, we have the entropy $$H(n)=\log (s+r)-2\log 2\frac{s}{s+r},$$ where $1\leq s\leq r$. We take $r\geq 3$, because we want to take at least three prime numbers in the decomposition of $n$ into prime factors and at least one square. We take a prime number $p$,
with gcd$(p,p_i)=1$ for all $i\in\{1,...,r\}$. We will study the difference of entropies $H(np^2)-H(np)$. This is $$H(np^2)-H(np)=\log \frac{s+r+2}{s+r+1}-2\log 2\frac{r+1}{(s+r+1)(s+r+2)}.$$ 
Next, using the Mathlab software program for different values of $s$, we deduce the values of $r$ for which $H(np^2)-H(np)<0$. Thus, we obtained the following list: $s=1$ and $r\geq 3$; $s=2$ and $r\geq 6$; $s=3$ and $r\geq 9$; $s=4$ and $r\geq 11$; $s=5$ and $r\geq 14$; $s=6$ and $r\geq 16$; $s=7$ and $r\geq 19$; $s=8$ and $r\geq 21$; $s=9$ and $r\geq 24$; $s=10$ and $r\geq 27$.\\
A plot of the function $f(s,r)=\log \dfrac{s+r+2}{s+r+1}-2\log 2\dfrac{r+1}{(s+r+1)(s+r+2)}$, with $s,r\in[0,100]$ is given below. 
\begin{center}
\includegraphics[scale = .8]{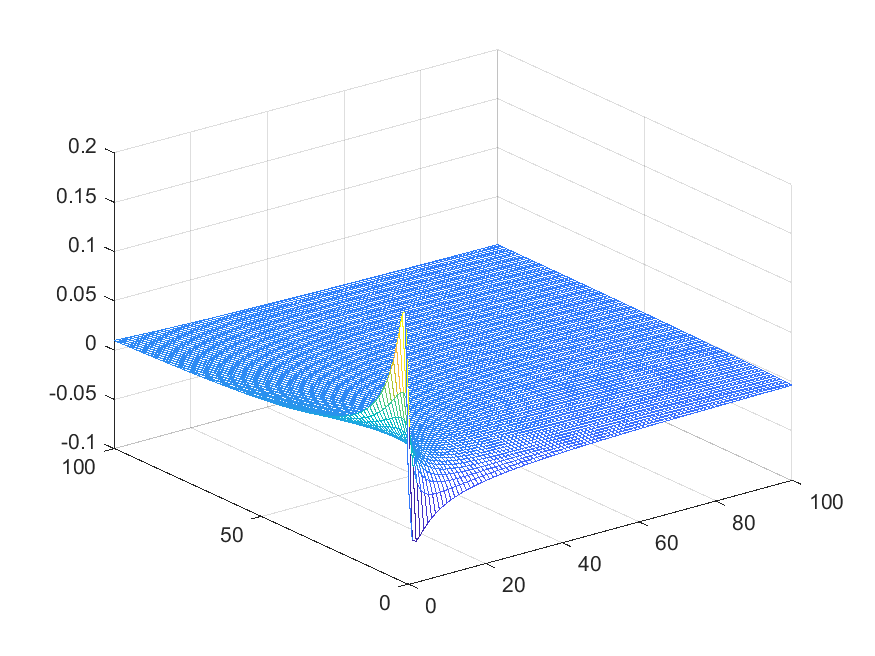}
\end{center}
For $r=s$ in decomposition of $n$ given above, we deduce that $H(np^2)-H(np)>0$.

We ask ourselves the problem of obtaining a general result.
\begin{proposition}
\label{new_1}
Let $n=p^{2}_{1}p^{2}_{2}...p_s^2p_{s+1}...p_{r}>1$ be an integer number, $1\leq s\leq r$ and $r\geq \frac{8s+5}{3}$. Then $H(np^2)-H(np)<0$, where $p$ is a prime number and $\text{gcd}(p,p_i)=1$ for all $i\in\{1,...,r\}$.
\end{proposition}
\begin{proof}
Using the Lagrange Theorem we deduce the following inequality: $$\log (x+1)-\log x<\frac{1}{x},$$ where $x>0$. Therefore, we have $$H(np^2)-H(np)=\log (s+r+2)-\log (s+r+1)-2\log 2\frac{r+1}{(s+r+1)(s+r+2)}$$ $$<\frac{1}{s+r+1}-2\log 2\frac{r+1}{(s+r+1)(s+r+2)}=\frac{s+r+2-2(\log 2)(r+1)}{(s+r+1)(s+r+2)}$$ $$\leq\frac{\frac{3r-5}{8}+r+2-2(\log 2)(r+1)}{(s+r+1)(s+r+2)}=\left(\frac{11}{8}-2\log 2\right)\frac{r+1}{(s+r+1)(s+r+2)}<0,$$
because $\frac{11}{8}-2\log 2=1.375-1.386...<0$. Consequently, we deduce the statement.
\end{proof}
\begin{remark}
With the assumptions from the statement of Proposition \ref{new_1}, we find the following inequality:
$H(np^2)-H(np)<0$, when we have $s=3k, r\ge 8k+2$ or $s=3k+1, r\ge 8k+5$ or $s=3k+2, r\ge8k+7$, with $k\geq 1$.
\end{remark}

If  $n=p^{\alpha_{1}}_{1}p^{\alpha_{2}}_{2}...p^{\alpha_{r}}_{r}$ and $m=p^{\alpha_{1}+\epsilon}_{1}p^{\alpha_{2}-\epsilon}_{2}...p^{\alpha_{r}}_{r}$, where $\epsilon\in\mathbb{N}$, $r,\alpha_{1}, \alpha_{2},...,\alpha_{r}\in\mathbb{N}^{*}$, $\alpha_2>\epsilon$ and  $p_{1}, p_{2},...,p_{r}$ are distinct prime positive integers. We remark that $\Omega(n)=\Omega(m)$.
It is easy to see that
\begin{equation}\label{eq_1}
H\left(m\right)-H\left(n\right)=\frac{1}{\Omega(n)}[\alpha_{1}\log\alpha_{1}+\alpha_{2}\log\alpha_{2}-(\alpha_{1}+\epsilon)\log(\alpha_{1}+\epsilon)-(\alpha_{2}-\epsilon)\log(\alpha_{2}-\epsilon)].
\end{equation}  
Therefore, our motivation is to study the difference in the entropies of the numbers $n=p^{\alpha}q^{\beta}$ and $m=p^{\alpha+\epsilon}q^{\beta-\epsilon}$, where $\epsilon\in\mathbb{N}$, $\alpha, \beta\in\mathbb{N}^{*}$, $\beta>\epsilon$ and  $p,q$ are distinct prime positive integers.
\begin{lemma}\label{Lemma 2.2}
Let $\alpha, \beta$ be two real numbers strictly positive. Then we have the inequality
\begin{equation}\label{ineq_1}
\frac{\alpha\log\alpha+\beta\log\beta}{\alpha+\beta}\geq\log\frac{\alpha+\beta}{2}.
\end{equation}
\end{lemma}
\begin{proof}
We consider the function $f:(0,\infty)\to\mathbb{R}$ defined by $f(x)=\alpha\log\alpha+x\log x-(x+\alpha)\log\frac{x+\alpha}{2}$. But, since $\frac{df}{dx}=\log\frac{2x}{\alpha+x}=0$, then $x=\alpha$. Since the function $f$ is decreasing on the interval $(0,\alpha]$ and increasing on $[\alpha,\infty)$, then $f(x)\geq f(\alpha)=0.$
\end{proof}
\begin{proposition}
\label{twodotone_1}
Let m,n be two numbers such that $n=p^{\alpha}q^{\beta}$ and $m=p^{\alpha+\epsilon}q^{\beta-\epsilon}$, with $\epsilon\in\mathbb{N}$, $\alpha, \beta\in\mathbb{N}^{*}$, $\beta>\epsilon$ and  $p,q$ are distinct prime positive integers. Then the inequality holds
\begin{equation}\label{ineq_1.1_0}
H\left(m\right)-H\left(n\right)\leq\frac{\alpha\log\alpha+\beta\log\beta}{\alpha+\beta}-\log\frac{\alpha+\beta}{2}.
\end{equation}
Moreover, if $\frac{\beta-\alpha}{2}\geq\epsilon$, then we have
\begin{equation}\label{ineq_1.1}
0\leq H\left(m\right)-H\left(n\right)\leq\frac{\alpha\log\alpha+\beta\log\beta}{\alpha+\beta}-\log\frac{\alpha+\beta}{2}.
\end{equation}
\end{proposition}
\begin{proof}
Using the definition of the entropy of a natural number, from \eqref{eq_1} for $r=2$, we obtain the following equality:
\begin{equation*}\label{eq_1.2}
H\left(m\right)-H\left(n\right)=\frac{1}{\alpha+\beta}[\alpha\log\alpha+\beta\log\beta-(\alpha+\epsilon)\log(\alpha+\epsilon)-(\beta-\epsilon)\log(\beta-\epsilon)].
\end{equation*}
From inequality \eqref{ineq_1}, replacing $\alpha$ and $\beta$ by $\alpha+\epsilon$ and $\beta-\epsilon$, we deduce
\begin{equation*}
(\alpha+\epsilon)\log(\alpha+\epsilon)+(\beta-\epsilon)\log(\beta-\epsilon)\geq(\alpha+\beta)\log\frac{\alpha+\beta}{2}.
\end{equation*}
Consequently, if we apply this inequality in the above equality, then we have the first inequality of the statement.

If $\alpha=\beta$, then from inequality $\frac{\beta-\alpha}{2}\geq\epsilon$, we deduce $\epsilon=0$, so $H\left(m\right)-H\left(n\right)=0$. Let $\alpha<\beta$, this implies $\beta>\frac{\beta-\alpha}{2}\geq\epsilon$. We take the function $f:[0,\frac{\beta-\alpha}{2}]\to\mathbb{R}$ defined by $f(t)=\alpha\log\alpha+\beta\log \beta-(\alpha+t)\log(\alpha+t)-(\beta-t)\log(\beta-t)$. Since $\frac{df}{dt}=\log\frac{\beta-t}{\alpha+t}=0$, then $t=\frac{\beta-\alpha}{2}$. The function $f$ is increasing on the interval $[0,\frac{\beta-\alpha}{2}]$, then $f(\frac{\beta-\alpha}{2})\geq f(t)\geq f(0)=0.$ Therefore, using the above equality and inequality \eqref{ineq_1.1_0}, we deduce inequality \eqref{ineq_1.1}.
\end{proof}
\begin{proposition}
\label{twodotone_2}
Let m,n,u be three numbers such that $n=p^{\alpha}q^{\beta}$ and $m=p^{\alpha+\epsilon}q^{\beta-\epsilon}$, with $\epsilon\in\mathbb{N}$, $\alpha, \beta\in\mathbb{N}^{*}$, $\beta>\epsilon$ and $gcd(m,u)=1,gcd(n,u)=1$, $p,q$ are distinct prime positive integers. Then the following inequality holds:
\begin{equation}\label{ineq_1.2}
H\left(mu\right)-H\left(nu\right)=\frac{\alpha+\beta}{\alpha+\beta+\Omega(u)}\left(H\left(m\right)-H\left(n\right)\right).
\end{equation}
\end{proposition}
\begin{proof}
Using the relation \eqref{eq_1} with $\Omega(mu)=\Omega(nu)=\alpha+\beta+\Omega(u)$
and the first equation in the proof of Proposition \ref{twodotone_1}, we deduce the equality of the statement.
\end{proof}
\begin{remark}
With the assumptions from the statement of Proposition \ref{twodotone_2}, we find the following inequality:
$$H\left(mu\right)-H\left(nu\right)\leq H\left(m\right)-H\left(n\right).$$
\end{remark}
Next, we will prove some results regarding the divergence of two numbers.

\begin{proposition}
\label{twodotone_1.1}
Let m,n be two numbers such that $n=p^{\alpha}q^{\beta}$ and $m=p^{\alpha+\epsilon}q^{\beta-\epsilon}$, with $\epsilon\in\mathbb{N}$, $\alpha, \beta\in\mathbb{N}^{*}$, $\beta>\epsilon$ and  $p<q$ are prime positive integers. Then the following inequality holds:
\begin{equation}\label{ineq_1.3}
D\left(n||m\right)\geq 0.
\end{equation}
\end{proposition}
\begin{proof}
If $\epsilon=0$, then $n=m$, so we have $D\left(n||m\right)=0$. We take $\epsilon>0$. From the definition of the divergence of two positive integers $n,m$, we find the equality
\begin{equation}\label{eq_1.3}
D\left(n||m\right)=\frac{1}{\alpha+\beta}[\alpha\log\alpha+\beta\log\beta-\alpha\log(\alpha+\epsilon)-\beta\log(\beta-\epsilon)].
\end{equation}
We consider the function $f:[0,\beta)\to\mathbb{R}$ defined by $f(t)=\alpha\log\alpha+\beta\log \beta-\alpha\log(\alpha+t)-\beta\log(\beta-t)$. Since $\frac{df}{dt}=\frac{t(\alpha+\beta)}{(\alpha+t)(\beta-t)}\geq 0$, then the function $f$ is increasing, so $f(t)\geq f(0)=0.$ Therefore, using equality \eqref{eq_1.3}, we have inequality \eqref{ineq_1.3}.
\end{proof}
\begin{proposition}
\label{twodotone_2.1}
Let m,n,u be three numbers such that $n=p^{\alpha}q^{\beta}$ and $m=p^{\alpha+\epsilon}q^{\beta-\epsilon}$, with $\epsilon\in\mathbb{N}$, $\alpha, \beta\in\mathbb{N}^{*}$, $\beta>\epsilon$ and $gcd(m,u)=1,gcd(n,u)=1$, $p<q$ are prime positive integers. Then the following inequality holds:
\begin{equation}\label{ineq_1.4}
D\left(nu||mu\right)=\frac{\alpha+\beta}{\alpha+\beta+\Omega(u)}D\left(n||m\right).
\end{equation}
\end{proposition}
\begin{proof}
Using relations \eqref{1.4} and \eqref{eq_1.3}, we deduce the equality of the statement.
\end{proof}
\begin{remark}
With the assumptions from the statement of Proposition \ref{twodotone_2.1}, we find the following inequality:
$$D\left(nu||mu\right)\leq D\left(n||m\right).$$
\end{remark}

Another problem that we want to study further is the determination of $m$ and $n$ when $D\left(n||m\right)=0$ knowing that gcd$(n,m)=1$ and $\Omega(n)=\Omega(m)$.

The Kullback--Leibler distance between two positive integer numbers $n,m\geq 2$ with  factorizations $n=p^{\alpha_{1}}_{1}p^{\alpha_{2}}_{2}...p^{\alpha_{r}}_{r}$ and $m=q^{\beta_{1}}_{1}q^{\beta_{2}}_{2}...q^{\beta_{r}}_{r}$, where the prime factors are ordered in ascending order, $\omega(n)=\omega(m)$ and $\Omega(n)=\Omega(m)$, as follows
$$                                            
D(n||m)=-\frac{1}{\Omega(n)}\sum_{i=1}^{r}\alpha_i\cdot \log \frac{\beta_i}{\alpha_i}. 
$$
It is easy to see that for $\alpha_i=\beta_i$ for all $i\in\{1,...,r\}$, we have $D\left(n||m\right)=0$. Therefore, we have to solve the system of equations
$\left\{\begin{matrix}
 \sum_{i=1}^{r}\alpha_i=\sum_{i=1}^{r}\beta_i \\ 
 \sum_{i=1}^{r}\alpha_i\cdot \log \frac{\beta_i}{\alpha_i}=0,  
\end{matrix}\right.$
with $\alpha_i\neq\beta_i$ for all $i\in\{1,...,r\}$.

For $r=2$, this system becomes
\begin{equation}\label{sys}
\left\{\begin{matrix}
\alpha_1+\alpha_2=\beta_1+\beta_2 \\ 
 \alpha_1^{\alpha_1}\alpha_2^{\alpha_2}=\beta_1^{\alpha_1}\beta_2^{\alpha_2}\\
 \alpha_1\neq\beta_1.  
\end{matrix}\right.
\end{equation}
The condition $\alpha_2\neq\beta_2$ is easily deduced from the fact that $\alpha_1\neq\beta_1$.

In the above system if $\alpha_1=\beta_2$, then we deduce from first equation of the system that $\alpha_2=\beta_1$. Thus, the second equation becomes $\alpha_1^{\alpha_1}\alpha_2^{\alpha_2}=\alpha_2^{\alpha_1}\alpha_1^{\alpha_2}$. Therefore, we obtain $\alpha_1=\alpha_2=\beta_1=\beta_2$, which is a contadiction.

In system \eqref{sys}, if we take $\alpha_1=\alpha_2$, then we obtain from first equation of the system that $2\alpha_1=\beta_1+\beta_2$. Thus, the second equation becomes $\alpha_1^{2}=\beta_1\beta_2$, so $(\beta_1+\beta_2)^2=4\beta_1\beta_2$. Therefore, we obtain $\beta_1=\beta_2$, so $\alpha_1=\alpha_2=\beta_1=\beta_2$,  which is a contadiction. Consequently, we have  $\alpha_1\neq\alpha_2$.
\begin{remark} 
If we look at this system with $\alpha_1,\alpha_2,\beta_1\in{\mathbb{N}^*}$ and $\beta_2\in{\mathbb{Z}}$, the system \eqref{sys} has an infinity of solutions given by $\alpha_1=\alpha,\alpha_2=2\alpha,\beta_1=4\alpha,\beta_2=-\alpha$, where $\alpha\in{\mathbb{N}^*}$.
\end{remark} 
Next, using the Mathlab software program and Magma software program for values $1\leq\alpha_1,\alpha_2,\beta_1,\beta_2\leq 4000$ we did not find any solution for system \eqref{sys}.
This observation suggested the remark, the system 
\begin{equation}\label{sys_1}
\left\{\begin{matrix}
x+y=u+v \\ 
 x^{x}y^{y}=u^{x}v^{y}  
\end{matrix}\right.
\end{equation}
has no solution, where $x,y,u,v\in\mathbb{N}^*$ such that $x\neq u$.

The second equation of system \eqref{sys_1} becomes:
\begin{equation}\label{eq.1}
 x^{x}y^{y}=u^{x}\left(x+y-u\right)^{y},  
\end{equation}
where $x,y,u\in\mathbb{N}^*$ such that $x\neq u$.

Next we will show that this equation has no solutions even for real numbers.

\begin{lemma} \label{lemm}
Let two real numbers $x,y>0$ and $x\neq 1$. The following equation:
\begin{equation}\label{eqn_1}
 x^{x}y^{y}=(x+y-1)^{y}  
\end{equation}
has no solution.
\end{lemma}
\begin{proof}
If $y=1$, then we have $x^x=x$. It follows that $x=1$, which is false, so we find that $y\neq 1$.  We are still studying the case when $x=y$, with $x\neq 1$. Equation \eqref{eqn_1} becomes $x^{2x}=(2x-1)^x$, so, $x^2=2x-1$, which gives the solution $x=1$, which is a contradiction. Consequently, $x\neq y$.

Next, we will study the following cases:

I) For $1<y<x$ relation \eqref{eqn_1} becomes $(x+y-1)^{y}=x^{x}y^{y}>x^{y}y^{y}=(xy)^{y}$. It follows that $x+y-1>xy$, which is equivalent to $0>(x-1)(y-1)$, which is false.

II) For $1<x<y$, by logarithmization we get
$x\log x+y\log y=y\log (x+y-1)$, which prove that $x\log x=y[\log (x+y-1)-\log y]$. For $x$ fixed, using Lagrange's Theorem, there is $\theta\in(y,y+x-1)$ such that $x\log x=y\frac{x-1}{\theta}$. Making the limit for $y\to\infty$, we deduce 
\begin{equation}\label{eqn_1.1}
 x\log x=x-1,  
\end{equation}
with $x>1$. Since the function $g:(1,\infty)\to\mathbb{R}$ defined by $g(x)=x\log x-x+1$ is strictly increasing on $(1,\infty)$ we deduce that $x\log x>x-1$. Therefore, equation \eqref{eqn_1.1} has no solution, when $x>1$. 

III) For $0<x<y<1$ relation \eqref{eqn_1} becomes $(x+y-1)^{y}=x^{x}y^{y}>x^{y}y^{y}=(xy)^{y}$. We deduce that $x+y-1>xy$, which is equivalent to $0>(x-1)(y-1)$, which is false.

IV) For $0<y<x<1$, by logarithmization we get 
\begin{equation}\label{eqn_1.2}
 x\log x+y\log y=y\log (x+y-1).  
\end{equation}
For $y$ fixed, we consider the function $h_1:(y,1)\to\mathbb{R}$ defined by $h_1(x)=y\log (x+y-1)-x\log x-y\log y$ is strictly increasing on $(y,1)$, because $h_1'(x)=\frac{1-x}{x+y-1}-\log x>0$. It follows that $y\log (x+y-1)-x\log x-y\log y<0$. Therefore, equation \eqref{eqn_1.2} has no solution, when $y<x<1$.

V) For $0<x<1<y$, by logarithmization we obtain relation \eqref{eqn_1.2}. For $y$ fixed, we consider the function $h_2:(0,1)\to\mathbb{R}$ defined by $h_2(x)=y\log (x+y-1)-x\log x-y\log y$ is strictly increasing on $(0,1)$, because $h_2'(x)=\frac{1-x}{x+y-1}-\log x>0$. It follows that $y\log (x+y-1)-x\log x-y\log y<0$. Therefore, equation \eqref{eqn_1.2} has no solution, when $0<x<1$.

VI) For $0<y<1<x$, by logarithmization we obtain relation \eqref{eqn_1.2}. For $y$ fixed, we consider the function $h_3:(1,\infty)\to\mathbb{R}$ defined by $h_3(x)=y\log (x+y-1)-x\log x-y\log y$ is strictly decreasing on $(1,\infty)$, because $h_2'(x)=\frac{1-x}{x+y-1}-\log x<0$. It follows that $y\log (x+y-1)-x\log x-y\log y<0$. Therefore, equation \eqref{eqn_1.2} has no solution, when $x>1$.

Consequently, the equation of the statement has no solution, when $x,y>0$ and $x\neq 1$. 
\end{proof}
\begin{theorem}\label{teo}
Let three real numbers $x,y,u>0$ and $x\neq u$. The following equation has no solution:
\begin{equation*}
 x^{x}y^{y}=u^{x}(x+y-u)^{y}.  
\end{equation*}
\end{theorem}
\begin{proof}
By dividing by $u^{x+y}$ in the relation from the statement we get $\left(\frac{x}{u}\right)^{x}\left(\frac{y}{u}\right)^{y}=(\frac{x}{u}+\frac{y}{u}-1)^{y}$. It follows that $\left(\frac{x}{u}\right)^{\frac{x}{u}}\left(\frac{y}{u}\right)^{\frac{y}{u}}=(\frac{x}{u}+\frac{y}{u}-1)^{\frac{y}{u}}$.
If we make the notations $x_1=\frac{x}{u}, y_1=\frac{y}{u}$, then the previous equation becomes
$x_1^{x_1}y_1^{y_1}=(x_1+y_1-1)^{y_1}$, with $x_1\neq 1$. From Lemma \ref{lemm}, we prove that the equation of the statement has no solution, when $x,y,u>0$ and $x\neq u$.
\end{proof}
\begin{remark}\label{rem}
Using Theorem \ref{teo}, the system \eqref{sys_1} has no solution, when $x,y,u,v\in\mathbb{N}^*$ with $x\neq u$.
\end{remark}
\begin{theorem} \label{teo_2}
For two positive integer numbers $n,m\geq 2$ with  factorizations $n=p^{\alpha_{1}}_{1}p^{\alpha_{2}}_{2}$ and $m=q^{\beta_{1}}_{1}q^{\beta_{2}}_{2}$ and $\alpha_{1}+\alpha_{2}=\beta_1+\beta_2$, $p_1<p_2$, $q_1<q_2$, it follows that                                           
$D(n||m)=0$ if and only if $\alpha_1=\beta_1$ and $\alpha_2=\beta_2$.
 
\end{theorem}
\begin{proof}
If $\alpha_1=\beta_1$ and $\alpha_2=\beta_2$, then it easy to see that $D(n||m)=0$.

If $D(n||m)=0$, this we obtain
$$
\left\{\begin{matrix}
\alpha_1+\alpha_2=\beta_1+\beta_2 \\ 
 \alpha_1^{\alpha_1}\alpha_2^{\alpha_2}=\beta_1^{\alpha_1}\beta_2^{\alpha_2}  
\end{matrix}\right.
$$
If $\alpha_1\neq\beta_1$, then from Remark \ref{rem}, this system has no solution. Therefore, we find $\alpha_1=\beta_1$, which prove that $\alpha_2=\beta_2$, so, we have the statement.
\end{proof}
\begin{remark} From Proposition \ref{twodotone_1.1} and Theorem \ref{teo_2}, we deduce that $D\left(n||m\right)> 0$ for
 two positive integers numbers $m,n$ such that $n=p^{\alpha}q^{\beta}$ and $m=p^{\alpha+\epsilon}q^{\beta-\epsilon}$, with $\alpha, \beta, \epsilon\in\mathbb{N}^{*}$, $\beta>\epsilon$ and  $p<q$ are prime positive integers. 
\end{remark}

\section{Some inequalities involving the entropy of an ideal of a ring of algebraic integers and the divergence of two ideals of a ring of algebraic integers}
\label{s:AlgebraicNumberField}

\noindent 
Now let \(K\) be an algebraic number field.
In \cite{MS1}, the authors obtained the following results about the entropy of an ideal or about the divergence of two ideals.

\begin{proposition}
\label{onedotfive}
 Let $K$ be an algebraic number field and let $I\neq(1)$ be an ideal of the ring $\mathcal{O}_{K}.$ Let  $ \omega\left(I\right)$  be  the number of distinct prime divisors of the ideal $I.$ Then:
 \begin{equation}
0\leq H\left(I\right)\leq \log \: \omega\left(I\right).  \tag{1.7}
\end{equation}

\end{proposition}

\begin{remark}
\label{onedotsix}
Let $K$ be an algebraic number field and let $I,J\neq \left(1\right)$ be two ideals of the ring  $\mathcal{O}_{K}$,  uniquely decomposed as $I=P^{e_{1}}_{1}\cdot P^{e_{2}}_{2}\cdot\ldots\cdot P^{e_{g}}_{g}$ and $J=Q^{e^{'}_{1}}_{1}\cdot Q^{e^{'}_{2}}_{2}\cdot\ldots\cdot Q^{e^{'}_{g}}_{g}$,
with $e_{1}, e_{2},\ldots,e_{g}, e^{'}_{1}, e^{'}_{2},\ldots,e^{'}_{g}$ positive integers, $P_{1}, P_{2},\ldots,P_{g}$ distinct prime ideals of the ring $\mathcal{O}_{K}$ and $Q_{1}, Q_{2},\ldots,Q_{g}$ distinct prime ideals of the ring $\mathcal{O}_{K}.$
If $e_{i}=e^{'}_{i},$ for $i=1,\ldots,g$,  then $D\left(I||J\right)=D\left(J||I\right)=0$.
\end{remark}

Since the proof of Proposition 11 in \cite{MS1}
only refers to the proof of Theorem 2 in \cite{Minculete1},
we give an independent proof of Proposition \ref{onedotfive}:
\begin{proof}
Since the quotients $0<\frac{e_i}{\Omega(I)}\le 1$, $i=1,\ldots,g$, in the expression for the entropy of an ideal $I\ne (1)$ in Definition \ref{onedotthree}
form a probability distribution associated to $I$, the logarithms are $\log \: \frac{e_i}{\Omega(I)}\leq 0$, and thus the entropy
$H(I)= - \sum_{i=1}^g\,\frac{e_i}{\Omega(I)} \log \: \frac{e_i}{\Omega(I)}\geq 0$ is non-negative.
For the proof of the optimal upper bound $H(I)\leq \log \: \omega(I)$ we use Formula \eqref{1.5}
and the Jensen inequality $f\left(\frac{1}{g}\sum_{i=1}^g\,e_i\right)\le\frac{1}{g}\sum_{i=1}^g\,f(e_i)$
for the function $f:(0,\infty)\to\mathbb{R}$, $x\mapsto x\log x$, which is convex downwards, since $f^{\prime\prime}(x)=\frac{1}{x}>0$ for $x>0$.
We have
$\left(\frac{1}{g}\sum_{i=1}^g\,e_i\right)\log \:\left(\frac{1}{g}\sum_{i=1}^g\,e_i\right)\le\frac{1}{g}\sum_{i=1}^g\,e_i\log \: e_i$.
By multiplication with $g$, this inequality becomes
$$\Omega(I)\left(\log \:\left(\sum_{i=1}^g\,e_i\right)-\log \: g\right)\le\sum_{i=1}^g\,e_i\log \: e_i,$$
and division by $\Omega(I)$ finally yields
$$H(I)=\log \: \Omega(I) - \frac{1}{\Omega(I)}\cdot \sum_{i=1}^{g}\,e_i\cdot \log \: e_i\le \log \: g = \log \: \omega(I).\qedhere$$
\end{proof}

In the case $\omega(I)=1$ of a prime ideal power $I=P^\alpha$, the maximal and minimal entropy coincides, since trivially $H(I)=0=\log \: \omega(I)$.
We show that the maximal entropy of composite ideals $I=P^{e_1}_1\cdot P^{e_2}_2\cdot\ldots\cdot P^{e_g}_g$ with at least two prime ideal divisors,
$g=\omega(I)\ge 2$, attains its maximum $\log \: \omega(I)$
precisely for equal exponents $e_1=e_2=\ldots=e_g$.
This supplements the items (iii) and (iv) of Proposition \ref{onedotone}.

\begin{proposition}
\label{onedotseven}
Let $K$ be an algebraic number field and let $\mathcal{O}_{K}$ be its ring of algebraic integers. 
Let $J=P^{e_1}_1\cdot P^{e_2}_2\cdot\ldots\cdot P^{e_g}_g$ be an ideal of a ring $\mathcal{O}_K$ with $g=\omega(J)\ge 2$.
Then $H(J)=\log \: \omega(J)$ if and only if $e_1=e_2=\ldots=e_g$.
\end{proposition}

\begin{proof}
By Formula \eqref{1.5}, the entropy of $J$ is
$H(J)=\log \: \Omega(J) - \frac{1}{\Omega(J)}\cdot \sum_{i=1}^{g}\,e_i\cdot \log \: e_i$,
where $\log$ is the natural logarithm and $\Omega(J)=\sum_{i=1}^g\,e_i$.

Sufficiency ($\Longleftarrow$):
If $e_1=e_2=\ldots=e_g=:e$, then $\Omega(J)=\sum_{i=1}^g\,e=g\cdot e$ and
$$H(J)=\log \: (g\cdot e) - \frac{1}{g\cdot e}\cdot \sum_{i=1}^{g}\,e\cdot \log \: e
=\log \: g+\log \: e-\frac{g\cdot e\cdot \log \: e}{g\cdot e}=\log \: \omega(J).$$

Necessity ($\Longrightarrow$):
We consider the $g$-variate function
$$f:(1,\infty)^g\to\mathbb{R},\
(x_1,\ldots,x_g)\mapsto \log \: \left(\sum_{i=1}^g\,x_i\right)-\frac{\sum_{i=1}^g\,x_i\log \: x_i}{\sum_{i=1}^g\,x_i}.$$
Since $\frac{\partial}{\partial x_j}\left(\sum_{i=1}^g\,x_i\log \: x_i\right)=1\cdot\log \: x_j+x_j\frac{1}{x_j}$,
the first partial derivatives of $f$ are
\begin{equation*}
\begin{aligned}
\frac{\partial f}{\partial x_j} &= \frac{1}{\sum_{i=1}^g\,x_i}\cdot 1 -
\left(\frac{1}{\sum_{i=1}^g\,x_i}(\log x_j+1)+\frac{-1}{\left(\sum_{i=1}^g\,x_i\right)^2}\sum_{i=1}^g\,x_i\log \: x_i\right)\\
&= \frac{\sum_{i=1}^g\,x_i-\sum_{i=1}^g\,x_i\log \: x_j-\sum_{i=1}^g\,x_i+\sum_{i=1}^g\,x_i\log \: x_i}{\left(\sum_{i=1}^g\,x_i\right)^2}\\
&= \frac{\sum_{i=1}^g\,x_i\left(\log \: x_i-\log \: x_j\right)}{\left(\sum_{i=1}^g\,x_i\right)^2},\quad \text{ for } j=1,\ldots,g.
\end{aligned}
\end{equation*}
They certainly vanish, when all variables are equal, $x_1=x_2=\ldots=x_g$.
If not all variables $x_i$ are equal, let $x_j$ be the minimum of them.
Then $x_j<x_i$ for at least one $1\le i\le g$, and thus the difference $\log \: x_i-\log \: x_j>0 $
and the entire sum $\sum_{i=1}^g\,x_i\left(\log \: x_i-\log \: x_j\right)$ is positive.
Therefore, equality of all variables is mandatory for an extremum of the function $f$.
\end{proof}

Next, we generalize Proposition \ref{twodotone_1.0}, for ideals in rings of algebraic integers.
\begin{proposition}
\label{3.0} Let $K$ be an algebraic number field and let $\mathcal{O}_{K}$ be its ring of algebraic integers. Let $J\ne (1)$ be an ideal of the ring $\mathcal{O}_{K}.$ 
Then $H\left(J\right)=0$ if and only if $J=P^{\alpha}$, where $P$ is a prime ideal of the ring  $\mathcal{O}_{K}$ and $\alpha$ is a positive integer.
\end{proposition}
\begin{proof}
According to the fundamental theorem of Dedekind rings, 
an ideal $J\ne (1)$ has a unique representation $J=P_1^{e_1}P_2^{e_2}\cdots P_g^{e_g}$
with at least one prime ideal divisor, that is, $g\ge 1$,
distinct prime ideal factors $P_1,P_2,\ldots,P_g$ of the ring $\mathcal{O}_{K}$,
and non-zero exponents $e_i\ge 1$ for $i=1,\ldots,g$.
By Definition \ref{onedotthree}, the entropy of $J$ is given by
$H(J)=-\sum_{i=1}^g\,p(e_i)\cdot\log \: p(e_i)$,
where $\log$ is the natural logarithm and the $p(e_i)=\frac{e_i}{\Omega(J)}$
with $\Omega(J)=\sum_{i=1}^g\,e_i$
form a particular probability distribution associated to $J$. 

Sufficiency ($\Longleftarrow$):
If $J=P^\alpha$, then $g=1$, $P=P_1$, $\alpha=e_1$, $\Omega(J)=\alpha$, and $p(\alpha)=\frac{\alpha}{\Omega(J)}=1$,
whence $H(J)=-p(\alpha)\cdot\log \: p(\alpha)=-1\cdot\log \: 1=0$.

Necessity ($\Longrightarrow$):
Suppose that $H(J)=0$ for $J=P_1^{e_1}P_2^{e_2}\cdots P_g^{e_g}$.
If we had more than one prime ideal, that is, $g\ge 2$,
then $\Omega(J)=e_1+e_2+\ldots+e_g>e_i$, $p(e_i)=\frac{e_i}{\Omega(J)}<1$, and $\log \: p(e_i)<0$,
for each $i=1,\ldots,g$.
Consequently, the entropy $H(J)=-\sum_{i=1}^g\,p(e_i)\cdot\log \: p(e_i)$
would be a sum of at least two positive terms $p(e_i)\cdot(-\log \: p(e_i))>0$,
in contradiction to the assumption that $H(J)=0$.
Thus $g=1$ and $J=P^\alpha$ with $P=P_1$, $\alpha=e_1$.

We mention another way to show the necessity:
taking into account Formula \eqref{1.5}, we have:
$$
H\left(J\right)=0  \Leftrightarrow
\log \: \Omega\left(J\right) = \frac{1}{\Omega\left(J\right)}\cdot \sum_{i=1}^{g}\,e_{i}\cdot \log \: e_{i} \Leftrightarrow
\Omega\left(J\right) \cdot \log \: \Omega\left(J\right) =  \sum_{i=1}^{g}\, \log \left(e^{e_{i}}_{i} \right)
$$
\begin{equation}
\label{3.1}
\Leftrightarrow \left(e_{1} + e_{2}+ \cdots + e_{g}\right)^{e_{1}+ e_{2}+\cdots+ e_{g}} = e_{1}^{e_{1}}\cdot e_{2}^{e_{2}}\cdot \cdots \cdot e_{g}^{e_{g}}. \tag{3.1}
\end{equation}
We try to solve the Diophantine equation \eqref{3.1}.

Since $e_{1},e_{2},\ldots,e_{g}$ are positive integers, the following equation
$$
\left(e_{1} + e_{2}+ \cdots + e_{g}\right)^{e_{1}+ e_{2}+\cdots+ e_{g}} =
$$
$$
\left(e_{1} + e_{2}+ \cdots + e_{g}\right)^{e_{1}}\cdot \ldots \cdot\left(e_{1} + e_{2}+\ldots + e_{g}\right)^{e_{g}} = e_{1}^{e_{1}}\cdot e_{2}^{e_{2}}\cdot \ldots \cdot e_{g}^{e_{g}}
$$
is impossible for $g\ge 2$, since $e_{1}+ e_{2}+\cdots+ e_{g}>e_{i}$ for each $i=1,\ldots,g$.
Equality is achieved if and only if $g=1$ such that $e_{1}\geq 1$
and Formula \eqref{3.1} degenerates to the triviality $e_{1}^{e_{1}}=e_{1}^{e_{1}}$.
If we denote $e_{1}=\alpha$ and $P_1=P$, then we obtain that $J=P^{\alpha}$.
\end{proof}

We want to see if there is an analogue of Proposition \ref{twodotone_1} for ideals in certain rings of algebraic integers,
that is, we are looking for fields of algebraic numbers $K$ and two ideals $I$ and $J$ of the ring  $\mathcal{O}_{K}$
so that $I$ and $J$ are ideals with the same two prime divisors and $\Omega\left(I\right) = \Omega\left(J\right)$. \\

We are looking for such an example, when $K=\mathbb{Q}\left(\xi\right)$ is a cyclotomic field. It is known that the ring of algebraic integers of $K$ is $\mathbb{Z}\left[\xi\right]$. We denote by $U\left(\mathbb{Z}\left[\xi\right]\right)$ the set the set of invertible elements of the ring $\mathbb{Z}\left[\xi\right]$. \\

First, we recall some results about cyclotomic fields.
\begin{theorem}
\label{theorem3.1}
(\cite{Savin}, \cite{Washington}) 
 \textit{Let }$n$\textit{\ be a positive integer, }%
$n\geq 3$. \textit{\ Let} $\xi$\textit{\ be a primitive root of order} $n$ of the unity and let $\mathbb{Q}\left(\xi\right)$ be the $n$th cyclotomic field. If $p$\textit{\ is a prime
positive integer,} $p$ does not divide $n$ \textit{\ and }$f$%
\textit{\ is the smallest positive integer such that} $p^{f}\equiv 1$  (mod $n$),
\textit{then we have }$p\mathbb{Z}\left[\xi\right]=P_{1}P_{2}....P_{r},$ 
\textit{where } $r=\frac{\varphi \left( n\right) }{f},\varphi $\textit{\ is
the Euler's function and }$P_{j},\,j=1,...,r$\textit{\ are different prime
ideals in the ring }$\mathbb{Z}[\xi ].$
\end{theorem}

\begin{corollary}
\label{3.2}
(\cite{Washington}) 
\textit{Let} $\xi $ \textit{be a primitive root of order} $n$ of the unity, where $n$ is  a positive integer, $n\geq 3$.
 Let  $\mathbb{Q}\left(\xi\right)$ be the $n$th cyclotomic field. Let $p$ be a prime positive integer. Then  $p$ splits completely in the ring $\mathbb{Z}\left[\xi\right]$ 
if and only if $p\equiv 1$ (mod $n$).

\end{corollary}

\begin{corollary}
\label{3.3}
(\cite{ireland}) 
\textit{Let} $\xi $ \textit{be a primitive root of order} $n$ of the unity, where $n$ is  a positive integer, $n\geq 3$.
 Let  $\mathbb{Q}\left(\xi\right)$ be the $n$th cyclotomic field. Let $p$ be a prime positive integer and let $P$ be a prime ideal in $\mathbb{Z}\left[\xi_{n}\right]$  such that $P\cap \mathbb{Z}=p \mathbb{Z}.$ If $p$ is odd then $P$ is ramified if and only if $p|n.$ If $p=2$
then $P$ is ramified if and only if $4|n.$

\end{corollary}

\begin{proposition}
\label{3.4} 
(\cite{Savin}) 
Let $p$ be a prime positive integer ane let $\xi $ \textit{be a primitive root of order} $p$ of the unity. Let  $\mathbb{Q}\left(\xi\right)$ be the $p$th cyclotomic field.
Then, the following statements are true:
\begin{enumerate}
 \item
$1-\xi$ is a prime element of the ring  $\mathbb{Z}\left[\xi\right]$;
\item
$p=u\cdot \left(1- \xi\right)^{p-1},$ where $u\in U\left(\mathbb{Z}\left[\xi\right]\right)$. 
 \end{enumerate}
\end{proposition}
We find the following example: let $\xi_{5}$ be a primitive root of order $5$ of the unity and let $K=\mathbb{Q}\left(\xi_{5}\right)$ be the $5$th cyclotomic field. It is known that the ring of algebraic integers of the field $K$, $\mathbb{Z}\left[\xi_{5}\right]$ is a principal domain.
We denote by Spec($\mathbb{Z}\left[\xi_{5}\right]$) the set of prime ideals of the ring $\mathbb{Z}\left[\xi_{5}\right]$.
We consider the following ideals of this ring: 
$I=10 \mathbb{Z}\left[\xi_{5}\right]=2\mathbb{Z}\left[\xi_{5}\right]\cdot 5\mathbb{Z}\left[\xi_{5}\right], \;  J=16 \left(1- \xi\right) \mathbb{Z}\left[\xi_{5}\right]=2^{4}\mathbb{Z}\left[\xi_{5}\right] \cdot  \left(1- \xi\right) \mathbb{Z}\left[\xi_{5}\right]$
and $ J^{'}=4 \left(1- \xi\right)^{3} \mathbb{Z}\left[\xi_{5}\right]$ and we want to decompose these ideals into products of prime ideals of the ring $\mathbb{Z}\left[\xi_{5}\right]$. It is known that  $\left(1- \xi\right) \mathbb{Z}\left[\xi_{5}\right]\in Spec(\mathbb{Z}\left[\xi_{5}\right]$).

Since $ord_{\left(\mathbb{Z}^{*}_{5}; \cdot\right)}\left(\overline{2}\right)=4$, applying Theorem \ref{theorem3.1} we have $r=\frac{\varphi \left( 5\right) }{4}=1$. It results that $2\mathbb{Z}\left[\xi_{5}\right]\in Spec(\mathbb{Z}\left[\xi_{5}\right]$). According to Proposition \ref{3.4}, $5=u\cdot \left(1- \xi\right)^{4},$ where $u\in U\left(\mathbb{Z}\left[\xi_{5}\right]\right),$ so, the ideal $5\mathbb{Z}\left[\xi_{5}\right]=\left( \left(1- \xi\right) \mathbb{Z}\left[\xi_{5}\right]\right)^{4}$.\\
Applying the Fundamental Theorem of Dedekind rings, it turns out that the ideals $I$ and $J$ decompose uniquely into the product of prime ideals in the ring $\mathbb{Z}\left[\xi_{n}\right]$ thus:
$$I=2\mathbb{Z}\left[\xi_{5}\right]\cdot \left( \left(1- \xi\right) \mathbb{Z}\left[\xi_{5}\right]\right)^{4}, \;  J=\left(2\mathbb{Z}\left[\xi_{5}\right]\right)^{4}\cdot  \left(1- \xi\right) \mathbb{Z}\left[\xi_{5}\right] $$
$$ \text{and}\;  J^{'}= \left(2\mathbb{Z}\left[\xi_{5}\right]\right)^{2} \cdot  \left(\left(1- \xi\right)\mathbb{Z}\left[\xi_{5}\right]\right)^{3}.$$

Considering $\epsilon=3$, we can write $J=\left(2\mathbb{Z}\left[\xi_{5}\right]\right)^{1+\epsilon}\cdot  \left(\left(1- \xi\right) \mathbb{Z}\left[\xi_{5}\right]\right)^{4-\epsilon},$ that is $\Omega\left(I\right) = \Omega\left(J\right)=5$. So, 
applying formula (1.5), we obtain $H\left(J\right) - H\left(I\right)=0.$\\
For $I$ and $J^{'},$ applying formula (1.5) it results that $H\left(I\right)= \log 5 - \frac{1}{5}\log 256 $ and $H\left(J^{'}\right)= \log 5 - \frac{1}{5}\log 108$.\\

We remark that $H\left(I\right)$ and  $H\left(J^{'}\right)$ satisfy the inequality in  Proposition \ref{twodotone_1}  (for $\epsilon=1$ $\alpha=1,$ $\beta =4$), that is
$$H\left(J^{'}\right) - H\left(I\right) = \frac{1}{5}\log\left( \frac{64}{27}\right) \leq \frac{1}{5}\log\left( \frac{8192}{3125}\right)= \frac{\alpha \log\alpha + \beta log\beta}{\alpha +  \beta} - \log\left(\frac{\alpha+\beta}{2}\right).$$
The result from the previous example (with $\epsilon=1$) can be generalized as follows:
\begin{proposition}
\label{3.5}
\textit{Let} $\xi_{5}$ \textit{be a primitive root of order} $5$ of the unity and let $K=\mathbb{Q}\left(\xi_{5}\right)$ be the $5$th cyclotomic field. Let $r$ be a positive integer, let $p,  p_{1},\ldots, p_{r}$ be distinct prime positive integers, $p\equiv 2$ or $3$ (mod $5$),
$p_{1}\equiv p_{2}\equiv \ldots \equiv p_{r}\equiv 1$ (mod $5$) and let the ideals $I_{1}=5p\mathbb{Z}\left[\xi_{5}\right],$ 
$J_{1}=\left(1- \xi\right)^{3} \cdot p^{2}\cdot \mathbb{Z}\left[\xi_{5}\right],$  $I_{2}=5p\cdot p_{1}p_{2}\cdot \ldots \cdot p_{r} \mathbb{Z}\left[\xi_{5}\right],$ $J_{2}=\left(1- \xi\right)^{3} \cdot p^{2}\cdot p_{1}p_{2}\cdot \ldots \cdot p_{r}\mathbb{Z}\left[\xi_{5}\right].$  
Then, the following statements hold:
 \begin{enumerate}
 \item
$0\leq H\left(J_{1}\right) - H\left(I_{1}\right) < 0.193$;
 \item
$0\leq H\left(J_{2}\right) - H\left(I_{2}\right) < 0.046$.
 \end{enumerate}
\end{proposition}
\begin{proof}
(i) Since $p\equiv 2$ or $3$ (mod $5$), it immediately follows that  $ord_{\left(\mathbb{Z}^{*}_{5}; \cdot\right)}\left(\overline{p}\right)=4$ and applying Theorem \ref{theorem3.1} it results that $p\mathbb{Z}\left[\xi_{5}\right]\in Spec(\mathbb{Z}\left[\xi_{5}\right]$).
According to Proposition \ref{3.4}, $1-\xi$ is a prime element of the ring  $\mathbb{Z}\left[\xi_{5}\right]$ and $5$ is totally ramified in $\mathbb{Z}\left[\xi_{5}\right]$, therefore, the ideals $I_{1}$ and $J_{1}$ decompose uniquely into  in the product of prime ideals of the ring $\mathbb{Z}\left[\xi_{5}\right]$ 
thus:
$$I_{1}=p\mathbb{Z}\left[\xi_{5}\right]\cdot \left( \left(1- \xi\right) \mathbb{Z}\left[\xi_{5}\right]\right)^{4}\; \text{and} \;  J_{1}=\left(p\mathbb{Z}\left[\xi_{5}\right]\right)^{2}\cdot  \left(\left(1- \xi\right) \mathbb{Z}\left[\xi_{5}\right]\right)^{3}.$$
Similar to the previous example, we obtain $0\leq H\left(J_{1}\right) - H\left(I_{1}\right) \leq \frac{1}{5}\log\left( \frac{8192}{3125}\right)=0.1927... $.\\
(ii) Since $p_{i}\equiv 1$ (mod $5$) $\left(\forall\right)$, $i=\overline{1,5}$, applying Corollary \ref{3.2},  $p_{i}$ split completely in the ring $\mathbb{Z}\left[\xi_{5}\right]$,  $i=\overline{1,5}$. So, for each $i=\overline{1,5}$, the ideal  $p_{i}\mathbb{Z}\left[\xi_{5}\right]$ decomposes uniquely into the product of prime ideals of the ring $\mathbb{Z}\left[\xi_{5}\right]$ thus:
$$p_{i}\mathbb{Z}\left[\xi_{5}\right]=P_{i1}\cdot P_{i2}\cdot P_{i3} \cdot P_{i4}, \; \text{where}\; P_{ij}\in Spec(\mathbb{Z}\left[\xi_{5}\right])\; \left(\forall\right), j=\overline{1,4}. $$
Taking into account this and i), it turns out that the ideals $I_{2}$ and $J_{2}$ decompose uniquely into  in the product of prime ideals of the ring $\mathbb{Z}\left[\xi_{5}\right]$ 
thus:
$$I_{2}=p\mathbb{Z}\left[\xi_{5}\right]\cdot \left( \left(1- \xi\right) \mathbb{Z}\left[\xi_{5}\right]\right)^{4}\cdot P_{11}\cdot P_{12}\cdot P_{13} \cdot P_{14} P_{21}\cdot P_{22}\cdot P_{23} \cdot P_{24} \dots P_{41}\cdot P_{42}\cdot P_{43} \cdot P_{44}$$
and
$$ J_{2}=\left( p\mathbb{Z}\left[\xi_{5}\right]\right)^{2}\cdot \left( \left(1- \xi\right) \mathbb{Z}\left[\xi_{5}\right]\right)^{3}\cdot P_{11}\cdot P_{12}\cdot P_{13} \cdot P_{14} P_{21}\cdot P_{22}\cdot P_{23} \cdot P_{24} \dots P_{41}\cdot P_{42}\cdot P_{43} \cdot P_{44}.$$
Applying formula (1.5) we have $ H\left(I_{2}\right) =\log\left(21\right)- \frac{4\log 4}{21}$ and $ H\left(J_{2}\right) =\log\left(21\right)- \frac{2\log 2+3\log 3}{21}$.
So, we obtain $0\leq H\left(J_{2}\right) - H\left(I_{2}\right)= \frac{1}{21}\cdot \log\left(\frac{64}{27}\right) \leq \frac{1}{21}\log\left( \frac{8192}{3125}\right)=0.0458...$.

\end{proof}

Proposition \ref{3.5} can be generalized as follows:

\begin{proposition}
\label{3.6}
\textit{Let} $q$ be a prime positive integer, $q\geq 5$, let $\xi$ \textit{be a primitive root of order} $q$ of the unity and let $K=\mathbb{Q}\left(\xi\right)$ be the $q$th cyclotomic field. Let $r$ be a positive integer, let $p,  p_{1},\ldots, p_{r}$ be distinct prime positive integers, $\overline{p}=\left(\mathbb{Z}^{*}_{q},\cdot\right)$ and $ord_{\left(\mathbb{Z}^{*}_{q},\cdot\right)}\left(\overline{p_{i}}\right)\neq q-1,$ $\left(\forall\right) i=\overline{1,r}.$ Let the ideals $I_{1}=qp\mathbb{Z}\left[\xi\right],$ 
$J_{1}=\left(1- \xi\right)^{q-2} \cdot p^{2}\cdot \mathbb{Z}\left[\xi\right],$  $I_{2}=qp\cdot p_{1}p_{2}\cdot \ldots \cdot p_{r} \mathbb{Z}\left[\xi\right],$ $J_{2}=\left(1- \xi\right)^{q-2} \cdot p^{2}\cdot p_{1}p_{2}\cdot \ldots \cdot p_{r}\mathbb{Z}\left[\xi\right].$  
Then, the following statements hold:
 \begin{enumerate}
 \item
$0\leq H\left(J_{1}\right) - H\left(I_{1}\right) \leq \frac{\left(q-1\right)\cdot \log\left(q-1\right)}{q} - \log \frac{q}{2}$;
 \item
$0\leq H\left(J_{2}\right) - H\left(I_{2}\right) \leq \frac{\left(q-1\right)\cdot \log\left(q-1\right)}{q} - \log \frac{ q }{2}$.

 \end{enumerate}
\end{proposition}
\begin{proof}

(i) Since $\overline{p}=\left(\mathbb{Z}^{*}_{q},\cdot\right)$, it immediately follows that  $ord_{\left(\mathbb{Z}^{*}_{q}; \cdot\right)}\left(\overline{p}\right)=q-1$. According to Theorem \ref{theorem3.1} it follows that $p\mathbb{Z}\left[\xi\right]\in Spec(\mathbb{Z}\left[\xi\right]$).
According to Proposition \ref{3.4}, $1-\xi$ is a prime element of the ring  $\mathbb{Z}\left[\xi\right]$ and $q$ is totally ramified in $\mathbb{Z}\left[\xi\right]$, therefore, the ideals $I_{1}$ and $J_{1}$ decompose uniquely into  in the product of prime ideals of the ring $\mathbb{Z}\left[\xi\right]$ 
thus:
$$I_{1}=p\mathbb{Z}\left[\xi\right]\cdot \left( \left(1- \xi\right) \mathbb{Z}\left[\xi\right]\right)^{q-1}\; \text{and} \;  J_{1}=\left(p\mathbb{Z}\left[\xi\right]\right)^{2}\cdot  \left(\left(1- \xi\right) \mathbb{Z}\left[\xi\right]\right)^{q-2}.$$
Applying (1.5 )we have

$$ H\left(J_{1}\right) - H\left(I_{1}\right)= \frac{q-1}{q}\cdot \log \left(q-1\right) -\frac{2\log 2 +  \left(q-2\right)\cdot \log \left(q-2\right) }{q} $$

From here, it follows that
$$0\leq H\left(J_{1}\right) - H\left(I_{1}\right) \leq \frac{q-1}{q}\cdot \log \left(q-1\right) - \log \frac{q}{2},$$ which is true from Lemma \ref{Lemma 2.2}.

(ii) Since $ord_{\left(\mathbb{Z}^{*}_{q},\cdot\right)}\left(\overline{p_{i}}\right)\neq q-1,$ $i=\overline{1,r}$, applying Theorem \ref{theorem3.1},  $p_{i}$ split in the ring $\mathbb{Z}\left[\xi\right]$,  $i=\overline{1,r}$. So, for each $i=\overline{1,r}$, the ideal  $p_{i}\mathbb{Z}\left[\xi\right]$ decomposes uniquely into the product of prime ideals of the ring $\mathbb{Z}\left[\xi\right]$ thus:
$$p_{i}\mathbb{Z}\left[\xi\right]=P_{i1}\cdot P_{i2}\cdot \ldots\cdot P_{is_{i}}, \; \text{where}\; P_{ij}\in Spec(\mathbb{Z}\left[\xi\right]), \; \left(\forall\right)\;  i=\overline{1, r},   \; \left(\forall\right)\;  j=\overline{1, s_{i}}, $$
where  $s_{i}= \frac{q-1}{f_{i}}, \;  f_{i}= ord_{\left(\mathbb{Z}^{*}_{q},\cdot\right)}\left(\overline{p_{i}}\right)$ and $P_{ij}$, $i=1,...,r,$ $j=\overline{1, s_{i}}$ are different prime ideals in the ring  $\mathbb{Z}\left[\xi\right]$.
Taking into account this and i), it turns out that the ideals $I_{2}$ and $J_{2}$ decompose uniquely into  in the product of prime ideals of the ring $\mathbb{Z}\left[\xi\right]$ 
thus:
$$I_{2}=p\mathbb{Z}\left[\xi\right]\cdot \left( \left(1- \xi\right) \mathbb{Z}\left[\xi_{5}\right]\right)^{q-1}\cdot P_{11}\cdot  \ldots \cdot P_{1s_{1}} \cdot  \ldots \cdot  P_{r1}\cdot\ldots \cdot P_{rs_{r}}$$
and
$$ J_{2}=\left( p\mathbb{Z}\left[\xi\right]\right)^{2}\cdot \left( \left(1- \xi\right) \mathbb{Z}\left[\xi\right]\right)^{q-2}\cdot  \cdot P_{11}\cdot  \ldots \cdot P_{1s_{1}} \cdot  \ldots \cdot  P_{r1}\cdot\ldots \cdot P_{rs_{r}}.$$
Applying formula (1.5) we have

$$ H\left(I_{2}\right) =\log\left(q+ s_{1}+\ldots+ s_{r}\right) - \frac{q-1}{q+ s_{1}+\ldots+ s_{r}}\cdot \log \left(q-1\right)$$

 and 
$$ H\left(J_{2}\right) =\log\left(q+ s_{1}+\ldots+ s_{r}\right)- \frac{\left(q-2\right)\cdot \log\left(q-2\right)+2\log 2}{q+ s_{1}+\ldots+ s_{r}}$$

So, we obtain 
$$0\leq H\left(J_{2}\right) - H\left(I_{2}\right)=\frac{q-1}{q+ s_{1}+\ldots+ s_{r}}\cdot \log \left(q-1\right) -  \frac{\left(q-2\right)\cdot \log\left(q-2\right)+2\log 2}{q+ s_{1}+\ldots+ s_{r}}.$$
But $\Omega\left(I_{2}\right) = \Omega\left(J_{2}\right) = q+ s_{1}+\ldots+ s_{r}$.
From here, it follows that
$$0\leq H\left(J_{2}\right) - H\left(I_{2}\right)= \frac{q}{q+ s_{1}+\ldots+ s_{r}}\left(\frac{q-1}{q} \log \left(q-1\right) -\frac{2\log 2 +  \left(q-2\right) \log \left(q-2\right) }{q}\right)$$
$$ =\frac{q}{q+ s_{1}+\ldots+ s_{r}}\left( H\left(J_{1}\right) - H\left(I_{1}\right)\right). $$
Applying (i), we obtain that
$$0\leq H\left(J_{2}\right) - H\left(I_{2}\right) \leq \frac{\left(q-1\right)\cdot \log\left(q-1\right)}{q} - \log \frac{ q }{2}.$$
\end{proof} 

Proposition \ref{3.5} (i)  and Proposition \ref{3.6}(i) confirm the fact that the inequality in Proposition \ref{twodotone_1} also works for the entropy of the ideals of a ring of algebraic integers.

\begin{proposition}
\label{3.7}
Let $K$ be an algebraic number field and let $\mathcal{O}_{K}$ be its ring of algebraic integers. Let  $I$ and $J$ be two ideals of the ring $\mathcal{O}_{K}$ such that $I=P^{\alpha}_{1}\cdot P^{\beta}_{2}$ and
$J=P^{\alpha+\epsilon}_{1}\cdot P^{\beta-\epsilon}_{2},$ where $P_{1}$, $P_{2}$ are distinct prime ideals of the ring $\mathcal{O}_{K}$  and $\epsilon\in\mathbb{N}$, $\alpha, \beta\in\mathbb{N}^{*}$, $\frac{\beta-\alpha}{2}\geq\epsilon$. Then the following inequality holds:
$$ 0\leq H\left(J\right)-H\left(I\right)\leq\frac{\alpha\log\alpha+\beta\log\beta}{\alpha+\beta}-\log\frac{\alpha+\beta}{2}.$$
\end{proposition}

\begin{proof} 
The proof is similar to the proof of the Proposition  \ref{twodotone_1} .
\end{proof} 

We asked ourselves if there are rings of algebraic integers, in which there are many ideal pairs whose divergence is equal to $0$.\\
Let a cubic field $K=\mathbb{Q}\left(\theta\right)$ where $(\theta$ is a root of an irreducible polynomial of the type  $f=X^{3} - aX + b$$\in\mathbb{Z}\left[X\right]$. 
In \cite{LN}, P. Llorente and E. Nart made a complete classification of how any prime integer $p$ decomposes into the product of primes in the ring of algebraic integers of the cubic field $K$.\\
Let $\Delta= 4a^{3}- 27b^{2}.$ If $m$$\in$$\mathbb{Z},$ we denote by $v_{p}\left(m\right)$ the greatest power $k$ with the property $p^{k}|m$. Let $s_{p}=\frac{\Delta}{p^{v_{p}\left(\Delta\right)}}.$
\begin{center}

\end{center}

\begin{proposition}
\label{3.8}

(a part of Theorem 1 from  \cite{LN}).  Let a cubic field $K=\mathbb{Q}\left(\theta\right)$ and let  $f=X^{3} - aX + b$$\in\mathbb{Z}\left[X\right]$ be the minimal polynomial of $\theta$. Let $p$ be a prime integer, $p\geq 5.$ Let $\mathcal{O}_{K}$ be the ring of algebraic integers
of the field $K.$
Then, the following statements are true:
\begin{enumerate} 
 \item
if $p|a,$ $p|b$ and $1=v_{p}\left(a\right)< v_{p}\left(b\right),$ then the ideal $p\mathcal{O}_{K}=$$P_{1}\cdot P^{2}_{2},$ where $P_{1}$ and $P_{2}$ are distinct prime ideals of the ring $\mathcal{O}_{K}$;
\item
if $p$ does not  divide $ab$ and $s_{p}$ is odd, then the ideal $p\mathcal{O}_{K}=$$P_{1}\cdot P^{2}_{2},$ where $P_{1}$ and $P_{2}$ are distinct prime ideals of the ring $\mathcal{O}_{K}$.

 \end{enumerate}
Moreover, these are the only cases when a prime integer $p\geq 5$ has the decomposition  $p\mathcal{O}_{K}=$$P_{1}\cdot P^{2}_{2}$ in the ring $\mathcal{O}_{K}$, where $P_{1}$ and $P_{2}$ are distinct prime ideals of the ring $\mathcal{O}_{K}$.

\end{proposition}

Using this Proposition, we obtain we quickly obtain the following result.

\begin{proposition}
\label{3.9}
Let a cubic field $K=\mathbb{Q}\left(\theta\right)$ and let  $f=X^{3} - aX + b$$\in\mathbb{Z}\left[X\right]$ be the minimal polynomial of $\theta$. Let $\mathcal{O}_{K}$ be the ring of algebraic integers
of the field $K.$ Let $p$ and $q$ be two distinct prime integers, $p\geq 5,$ $q\geq 5.$  If $p$ and $q$ satisfy the conditions of hypothesis i) or the conditions of hypothesis ii) of the previous Proposition, then the following statements are true:\\
a) the entropies of the ideals $p\mathcal{O}_{K}$ and $q\mathcal{O}_{K}$ are equal;\\
b) the divergence $D(p\mathcal{O}_{K}||q\mathcal{O}_{K})=0$.

\end{proposition}

\begin{proof} 
a) The proof follows immediately, using Proposition  \ref{3.8}  and  formula (1.5).\\
b) The proof follows immediately, using Proposition  \ref{3.8}  and formula (1.6).
\end{proof}

\section{The entropy of probability measures for \(3\)-class field tower groups of imaginary quadratic number fields}
\label{s:ClassFieldTheory}

\noindent
According to Koch and Venkov
\cite{KoVe1975},
the Galois group \(\mathrm{Gal}(K^\infty/K)\) of the maximal unramified pro-\(3\)-extension \(K^\infty\)
of an imaginary quadratic field \(K=\mathbb{Q}(\sqrt{d})\) with fundamental discriminant \(d<0\)
must be a \textit{Schur \(\sigma\)-group} \(G\)
with balanced presentation, expressed by the coincidence of the
\textit{generator rank} \(d_1(G)=\mathrm{dim}_{\mathbb{F}_3}H^1(G,\mathbb{F}_3)\) and the
\textit{relation rank} \(d_2(G)=\mathrm{dim}_{\mathbb{F}_3}H^2(G,\mathbb{F}_3)\),
and with a \(\sigma\)-\textit{automorphism} \(\sigma\in\mathrm{Aut}(G)\),
acting by inversion \(x\mapsto x^{-1}\) on the cohomology groups
\(H^1(G,\mathbb{F}_3)\) and \(H^2(G,\mathbb{F}_3)\).
Denote by \(\mathrm{Cl}_3(K)=\mathrm{Syl}_3\mathrm{Cl}(K)\)
the \(3\)-class group of \(K\).
For \(3\)-class rank \(d_1(\mathrm{Cl}_3(K))=1\),
the \(3\)-class field tower \(K^\infty=K^1\) stops at the first stage with the Hilbert \(3\)-class field of \(K\),
and the group \(G=\mathrm{Gal}(K^1/K)\simeq\mathrm{Cl}_3(K)\) is non-trivial cyclic.
For \(d_1(\mathrm{Cl}_3(K))\ge 3\),
the tower \(K^\infty\) has infinitely many stages,
and \(G=\mathrm{Gal}(K^\infty/K)\) is a topological pro-\(3\)-group
\cite{KoVe1975}.
We study the intermediate situation with \(3\)-class rank
\(d_1(\mathrm{Cl}_3(K))=2\),
additionally assuming an elementary bicyclic \(3\)-class group
\(\mathrm{Cl}_3(K)=(\mathbb{Z}/3\mathbb{Z})^2\).
Consequently, we have to seek suitable Schur \(\sigma\)-groups \(G\)
on the descendant tree \(\mathcal{T}\) of the abelian root
\((\mathbb{Z}/3\mathbb{Z})^2=\langle 3^2,2\rangle\),
in the notation with order and identifier in angle brackets
of the SmallGroups database
\cite{BEO2002,BEO2005}.
We rigorously restrict the tree to
possible ancestors \(G\) of Schur \(\sigma\)-groups,
and thus we arrive at the pruned subtree \(\mathcal{T}_0\subset\mathcal{T}\)
which is drawn on page 657 of
\cite[\S\ 2.4]{BBH2017}.
In this tree diagram,
Boston, Bush and Hajir display the \textit{probability measure} \(p(G)\), defined in
\cite[Thm. 2.25, p. 653]{BBH2017},
of all relevant descendants \(G\) of the root \(\langle 3^2,2\rangle\)
in the range of orders \(3^3\le\mathrm{ord}(G)\le3^{12}\)
as rational fractions adjacent to the vertex which represents \(G\)
by a small full disc.
Schur \(\sigma\)-groups are surrounded additionally by a bigger contour circle.
From top to bottom,
the seven layers of vertices in the tree diagram have orders
\(3^3,3^5,3^6,3^8,3^9,3^{11},3^{12}\).
Since their descendant subtrees \(\mathcal{T}(R)\subset\mathcal{T}_0\)
are of eminent importance in class field theory,
we focus our attention on three non-abelian roots \(R=\langle 3^5,i\rangle\)
with \(i\in\lbrace 4,6,8\rbrace\) among the seven vertices with \(3\le i\le 9\).
Ascione et al.
\cite{AHL1977}
denote the five unique immediate \(\sigma\)-descendants by capital letters
\(B=\langle 3^6,40\rangle\), \(N=\langle 3^6,45\rangle\), \(Q=\langle 3^6,49\rangle\), \(U=\langle 3^6,54\rangle\), \(W=\langle 3^6,57\rangle\),
and call them \textit{non-CF groups} of \textit{second maximal nilpotency class} \(\mathrm{cl}\) or \textit{coclass} \(\mathrm{cc}=2\). 
For vertices \(G\) of the three infinite subtrees \(\mathcal{T}(R)\),
the probability measures \(p(G)\) obey rather amenable deterministic laws,
parametrized with indices \(n=0,1,2,\ldots\)
associated to the \textit{periodic structure} of these trees
\cite[Thm. 7.1, p. 167]{Ma2015a}. 

In terms of relative identifiers \(-\#s;i\) with step size \(s\),
which are given by the ANUPQ-package
\cite{GNO2006,GAP2025,MAGMA2025,MAGMA6561},
the periodicity is expressed by the formula
\begin{equation}
\label{eqn:TreeN}
G=G(n,j)=\langle 3^5,4\rangle(-\#1;1-\#2;j)^n, \quad n\ge 0, \quad j\in\lbrace 1,2\rbrace,
\end{equation}
for the infinite subtree \(\mathcal{T}(N)\) with finite branches,
where the vertex with \(j=1\) is extendible (capable),
and the vertex with \(j=2\) is a terminal Schur \(\sigma\)-group.
See the tree diagram in
\cite[Fig. 6, p. 154]{Ma2017},
where terminal descendants of step size \(s=1\) can be cancelled,
since they are not Schur \(\sigma\)-ancestors.
According to
\cite[Thm. 2.11(i), p. 645]{BBH2017},
the measures \(p(D)\) of the infinite collection of Schur \(\sigma\)-descendants \(D\)
of a vertex \(G\in\mathcal{T}(N)\) recursively sum up to the measure \(p(G)\),
and now we determine the entropy \(H(p)\) of this countable probability distribution \(p\).

\begin{theorem}
\label{thm:TreeN}
The sum of the absolute probability distribution \(p:\,S\to(0;1\rbrack\), \(G\mapsto p(G)\),
of the countable collection \(S\) of all Schur \(\sigma\)-groups \(G\) on the tree \(\mathcal{T}(N)\)
with root \(N=\langle 3^6,45\rangle\) is given by
\(\sum_{G\in S}\,p(G)=2^6\cdot 3^{-6}\approx 0,08779\).
It must be viewed with respect to all Schur \(\sigma\)-groups \(G\)
with order a power of the prime \(3\) and generator rank \(d_1(G)=2\).
The entropy of the normalized relative probability distribution
\(p_{\mathrm{rel}}:=\frac{p}{2^6\cdot 3^{-6}}\)
with \(\sum_{G\in S}\,p_{\mathrm{rel}}(G)=1\)
is given by (compare \eqref{eqn:EntropyU})
\begin{equation}
\label{eqn:EntropyN}
H(p_{\mathrm{rel}})=-\sum_{G\in S}\,p_{\mathrm{rel}}(G)\cdot\log(p_{\mathrm{rel}}(G))=\frac{3}{2}\log(3)-\log(2)\approx 0,9548.
\end{equation}
\end{theorem}

\begin{proof}
According to the formula in
\cite[Cor. 2.26, p. 655]{BBH2017},
where we insert the particular values \(p=3\) for \(3\)-groups and \(g=2\) for two generators,
the probability measure of a Schur \(\sigma\)-group \(G\) is given by
\begin{equation}
\label{eqn:Meas}
p(G)=\frac{y(G)^2}{\#\mathrm{Aut}(G)}\cdot 3^{2^2}\cdot\prod_{k=1}^2\,(1-\frac{1}{3^k})^2
\end{equation}
with a constant factor
\(3^4\cdot (1-\frac{1}{3})^2\cdot (1-\frac{1}{3^2})^2=3^4\cdot\frac{2^2}{3^2}\cdot\frac{8^2}{9^2}=\frac{2^8}{3^2}\)
and variable contributions by the number of fixed points \(y(G)\) of the \(\sigma\)-automorphism
and the order of the automorphism group \(\#\mathrm{Aut}(G)\).
For the tree \(\mathcal{T}(N)\),
all vertices share a common \textit{transfer kernel type} (TKT) H.4, \(\varkappa\sim (4111)\),
and we have the simple parametrized formation laws
\(y(G)^2=3^{2n+4}\) and \(\#\mathrm{Aut}(G)=2\cdot 3^{3n+9}\).
Together this yields
\(p(G)=\frac{2^8}{3^2}\cdot\frac{3^{2n+4}}{2\cdot 3^{3n+9}}=2^7\cdot 3^{-n-7}\)
with \(n\ge 0\).
Now we come to the total summation along the full infinite subtree \(\mathcal{T}(N)\), where we use geometric series: \\
\(\sum_{n=0}^\infty\,2^7 3^{-n-7}=\frac{128}{2187}\cdot\frac{3}{2}=\frac{64}{729}=\mathrm{meas}_4(N)\),
in the recursive sense of the class-\(c\) measure of the root \(N\) with nilpotency class \(c=\mathrm{cl}(N)=4\)
\cite[Thm. 2.13, p. 646]{BBH2017}.
Since the \textit{abelian quotient invariants} (AQI)
\((\lbrack 3,3\rbrack;\lbrack 3,3,3\rbrack^3,\lbrack 9,3\rbrack)\)
remain stable for all vertices of the tree \(\mathcal{T}(N)\),
the value \(\mathrm{meas}_4(N)=\frac{64}{729}\) precisely coincides with the measure of this IPAD in
\cite[Thm. 4.3(3), pp. 661--662]{BBH2017}.
For the entropy, we need the logarithms of the relative probabilities
\(p_{\mathrm{rel}}(G)=\frac{2^7\cdot 3^{-n-7}}{2^6\cdot 3^{-6}}=\frac{2}{3^{n+1}}\), that is,
\(\log(p_{\mathrm{rel}}(G))=\log(2)-(n+1)\log(3)\).
Therefore \\
\(H(p_{\mathrm{rel}})=-\sum_{G\in S}\,p_{\mathrm{rel}}(G)\cdot\log(p_{\mathrm{rel}}(G))
=\sum_{n=0}^\infty\,\frac{2}{3^{n+1}}\cdot\bigl((n+1)\log(3)-\log(2)\bigr)\)\\
\(=\frac{2}{3}\cdot\bigl(\log(3)\cdot\sum_{n=0}^\infty\,\frac{n}{3^n}+(\log(3)-\log(2))\cdot\sum_{n=0}^\infty\,\frac{1}{3^n}\bigr)\) \\
\(=\frac{2}{3}\cdot\bigl(\log(3)\cdot\frac{3}{4}+(\log(3)-\log(2))\cdot\frac{3}{2}\bigr)=\frac{3}{2}\log(3)-\log(2)\approx 0,9548\).
\end{proof}

\begin{remark}
\label{Rmk:TreeN}
The periodic structure of the tree \(\mathcal{T}(N)\) was also analyzed in
\cite[\S\ 6.2.2, pp. 299--304]{Ma2015b}
where a tree diagram
with Schur \(\sigma\)-groups \(S_0,S_1,S_2,S_3\)
is drawn in Figure 1 on page 302.
Bartholdi and Bush
\cite{BaBu2007}
have shown that the soluble length of the countable collection of Schur \(\sigma\)-groups in \(\mathcal{T}(N)\)
is unbounded, for instance \(\mathrm{sl}(S_i)=3\) for \(i=0,1,2\), but \(\mathrm{sl}(S_3)=4\).
\end{remark}

\noindent
The infinite subtrees \(\mathcal{T}(Q)\) and \(\mathcal{T}(U)\)
are isomorphic as digraphs.
The periodicity of \(\mathcal{T}(U)\) is of considerably higher complexity, expressed by the formula
\begin{equation}
\label{eqn:TreeU}
G=G(n,j)=\langle 3^5,8\rangle(-\#1;1-\#2;j)^n, \quad n\ge 0, \quad j\in\lbrace 1,\ldots,6\rbrace,
\end{equation}
where the vertices with \(j=1,2,3\) are extendible (capable),
one of them, \(j=1\), mainline with TKT c.21, \(\varkappa\sim (2034)\),
two of them, \(j=2,3\), with TKT G.16, \(\varkappa\sim (2134)\),
roots of \textit{infinite branches},
and the vertices with \(j=4,5,6\) are terminal Schur \(\sigma\)-groups,
two of them, \(j=4,6\), with TKT E.9, \(\varkappa\sim (2334)\sim (2434)\),
and one, \(j=5\), with TKT E.8, \(\varkappa\sim (2234)\).
The next theorem holds also for \(\mathcal{T}(Q)\).

\begin{theorem}
\label{thm:TreeU}
The sum of the absolute probability distribution \(p:\,S\to(0;1\rbrack\), \(G\mapsto p(G)\),
of the countable collection \(S\) of all Schur \(\sigma\)-groups \(G\) on the tree \(\mathcal{T}(U)\)
with root \(U=\langle 3^6,54\rangle\) is given by
\(\sum_{G\in S}\,p(G)=2^6\cdot 3^{-6}\approx 0,08779\).
It must be viewed with respect to all Schur \(\sigma\)-groups \(G\)
with order a power of the prime \(3\) and generator rank \(d_1(G)=2\).
The entropy of the normalized relative probability distribution
\(p_{\mathrm{rel}}:=\frac{p}{2^6\cdot 3^{-6}}\)
with \(\sum_{G\in S}\,p_{\mathrm{rel}}(G)=1\)
is given by (compare \eqref{eqn:EntropyN})
\begin{equation}
\label{eqn:EntropyU}
H(p_{\mathrm{rel}})=-\sum_{G\in S}\,p_{\mathrm{rel}}(G)\cdot\log(p_{\mathrm{rel}}(G))=\frac{33}{32}\log(3)-\frac{3}{8}\log(2)\approx 0,8730.
\end{equation}
\end{theorem}

\begin{proof}
Again we apply the formula
\eqref{eqn:Meas}.
However, for the tree \(\mathcal{T}(U)\), we have more complicated parametrized formation laws
with two different contributions,
\(y(G)^2=3^{2n+4}\) and \(\#\mathrm{Aut}(G)=2\cdot 3^{4n+10}\)
for the finite branches with TKT E.8, \(\varkappa\sim (2234)\), and TKT E.9, \(\varkappa\sim (2334)\sim (2434)\),
but \(y(G)^2=3^{2n+m+6}\) and \(\#\mathrm{Aut}(G)=2\cdot 3^{4n+3m+13}\)
for the infinite branches with TKT G.16, \(\varkappa\sim (2134)\).
Together this yields
\(p(G)=\frac{2^8}{3^2}\cdot\frac{3^{2n+4}}{2\cdot 3^{4n+10}}=2^7\cdot 3^{-2n-8}\) for TKT E.8, E.9, and
\(p(G)=\frac{2^8}{3^2}\cdot\frac{3^{2n+m+6}}{2\cdot 3^{4n+3m+13}}=2^7\cdot 3^{-2n-m-9}\) for TKT G.16,
with \(n\ge 0\), \(m\ge 0\).
Now we come to the total summation along the full infinite subtree \(\mathcal{T}(U)\),
including all infinite branches.
Again we use geometric series: \\
Firstly, \(3\cdot\sum_{n=0}^\infty\,2^7 3^{-2n-8}=\frac{128}{2187}\cdot\frac{9}{8}=\frac{16}{243}\),
for three Schur \(\sigma\)-groups with TKT E.8 and E.9 in each period.
Secondly, \(2\cdot\sum_{n=0}^\infty\,\sum_{m=0}^\infty\,2^7 3^{-2n-m-9}=\frac{256}{19683}\cdot\frac{3}{2}\cdot\frac{9}{8}=\frac{16}{729}\),
for all Schur \(\sigma\)-groups with TKT G.16 on two infinite branches arising in each period.
Together \(\frac{16}{243}+\frac{16}{729}=\frac{16}{729}\cdot (3+1)=\frac{64}{729}=\mathrm{meas}_4(U)\),
in the recursive sense of the class-\(c\) measure of the root \(U\) with nilpotency class \(c=\mathrm{cl}(U)=4\)
\cite[Thm. 2.13, p. 646]{BBH2017},
which coincides with the sum of parametrized IPAD measures in
\cite[Thm. 4.3(5,6), p. 662]{BBH2017}.
For the entropy, we need the logarithms of the relative probabilities.
Firstly,
\(p_{\mathrm{rel}}(G)=\frac{2^7\cdot 3^{-2n-8}}{2^6\cdot 3^{-6}}=\frac{2}{3^{2n+2}}\), that is,
\(\log(p_{\mathrm{rel}}(G))=\log(2)-(2n+2)\log(3)\), for TKT E.8 and E.9.
Secondly,
\(p_{\mathrm{rel}}(G)=\frac{2^7\cdot 3^{-2n-m-9}}{2^6\cdot 3^{-6}}=\frac{2}{3^{2n+m+3}}\), that is,
\(\log(p_{\mathrm{rel}}(G))=\log(2)-(2n+m+3)\log(3)\), for TKT G.16.
Therefore \\
\(H(p_{\mathrm{rel}})=-\sum_{G\in S}\,p_{\mathrm{rel}}(G)\cdot\log(p_{\mathrm{rel}}(G))
=\sum_{n=0}^\infty\,\frac{2}{3^{2n+2}}\cdot\bigl((2n+2)\log(3)-\log(2)\bigr)\) \\
\(+\sum_{n=0}^\infty\,\sum_{m=0}^\infty\,\frac{2}{3^{2n+m+3}}\cdot\bigl((2n+m+3)\log(3)-\log(2)\bigr)\) \\
\(=\frac{2}{9}\cdot\bigl(2\log(3)\cdot\sum_{n=0}^\infty\,\frac{n}{9^n}+(2\log(3)-\log(2))\cdot\sum_{n=0}^\infty\,\frac{1}{9^n}\bigr)\) \\
\(+\frac{2}{27}\cdot\bigl(2\log(3)\cdot\sum_{n=0}^\infty\,\sum_{m=0}^\infty\,\frac{n}{9^n}\cdot\frac{1}{3^m}+\log(3)\cdot\sum_{n=0}^\infty\,\sum_{m=0}^\infty\,\frac{1}{9^n}\cdot\frac{m}{3^m}\) \\
\(+(3\log(3)-\log(2))\cdot\sum_{n=0}^\infty\,\sum_{m=0}^\infty\,\frac{1}{9^n}\cdot\frac{1}{3^m}\bigr)\) \\
\(=\frac{2}{9}\cdot\bigl(2\log(3)\cdot\frac{9}{64}+(2\log(3)-\log(2))\cdot\frac{9}{8}\bigr)\) \\
\(+\frac{2}{27}\cdot\bigl(2\log(3)\cdot\frac{3}{2}\cdot\frac{9}{64}+\log(3)\cdot\frac{3}{4}\cdot\frac{9}{8}+(3\log(3)-\log(2))\cdot\frac{3}{2}\cdot\frac{9}{8}\bigr)\) \\
\(=\frac{1}{16}\log(3)+\frac{1}{4}(2\log(3)-\log(2))+\frac{1}{32}\log(3)+\frac{1}{16}\log(3)+\frac{1}{8}(3\log(3)-\log(2))\) \\
\(=\frac{2+16+1+2+12}{32}\log(3)-\frac{2+1}{8}\log(2)=\frac{33}{32}\log(3)-\frac{3}{8}\log(2)\approx 0,8730\).
\end{proof}

\begin{remark}
\label{Rmk:TreeU}
The periodic structure of the trees \(\mathcal{T}(Q)\) and \(\mathcal{T}(U)\) was also analyzed in
\cite[\S\ 6.2.2, pp. 184--193]{Ma2015a}
where tree diagrams
with Schur \(\sigma\)-groups \(\langle 3^8,i\rangle\), \(i\in\lbrace 616,617,618,620,622,624\rbrace\),
and others with bigger orders \(3^{11}\) and \(3^{14}\) 
are drawn in Figures 8 and 9 on pp. 188--189,
indicated with relative identifiers
\(Q-\#2;j\), \(j=4,5,6\), and \(U-\#2;j\), \(j=2,4,6\).
The soluble length of all these Schur \(\sigma\)-groups \(G\) with TKT E.6, E.14, E.8, E.9
is uniformly bounded by \(\mathrm{sl}(G)=3\).
Bush and Mayer
\cite{BuMa2015}
have shown that the two non-metabelian Schur \(\sigma\)-groups with TKT E.9 and \(i=620,624\)
disprove the erroneous claim by Scholz and Taussky
\cite[p. 41]{SoTa1934}
that \(\mathbb{Q}(\sqrt{-9748})\) has a metabelian \(3\)-class field tower with two stages.
The infinite branches with TKT G.16 were intentionally cancelled in Figure 9,
however, the corresponding infinite branches with TKT H.4 which are purged in Figure 8,
are discussed in context with Figure 4 in
\cite[pp. 101--102]{Ma2018}.
\end{remark}

\noindent
{\bf Conclusions:}
The entropy
\(H(p)=\sum_{v\in S}\,p(v)\cdot\log(p(v))\)
of a probability distribution \(p\) associated
with numbers in \S\ 2 and with ideals in \S\ 3,
where \(S\) is a \textit{finite set} of places,
turned out to take its \textit{maximum} \(\log(\lvert S\rvert)\)
for a \textit{Laplace distribution} with constant value \(\lvert S\rvert^{-1}\),
corresponding to maximal disorder,
and its \textit{minimum} \(0\)
for a \textit{Dirac-} (or Kronecker-)\textit{distribution} with sharp localization in a singleton \(\lvert S\rvert=1\),
corresponding to maximal order.
For a probability measure \(p\) on Schur \(\sigma\)-groups, however,
\S\ 4 shows that the extreme scenarios 
of the entropy with a \textit{countable set} \(S\) of tree vertices are disabled,
since Laplace-equidistribution and sharp Dirac-distribution do not exist.
Astonishingly,
the entropy \(0,8730\) of the complicated trees \(\mathcal{T}(Q)\) and \(\mathcal{T}(U)\)
is less than \(0,9548\) for the simple tree \(\mathcal{T}(N)\).


\noindent
{\bf Acknowledgments:} The authors are grateful to Professor Constantin Gheorghies (from Dunarea de Jos University of Gala\c{t}i) for interesting discussions related to this topic.


\end{document}